\documentclass{article}

\usepackage{amsmath,amssymb, enumerate}
\usepackage{graphicx, float}
\usepackage{enumitem, enumerate}
\usepackage{dsfont}
\usepackage[utf8]{inputenc}
\usepackage{tikz}
\usetikzlibrary{matrix,shapes.geometric,arrows,arrows.meta,positioning,chains}
\usepackage{physics}
\usepackage{amsthm}
\usepackage[]{algorithm2e}
\usepackage{tabularx}
\usepackage{fourier}
\usepackage{array}
\usepackage{booktabs}
\usepackage{tabularx}
\usepackage{makecell}

\newcommand{\comment}[1]{}
\newtheorem*{theorem*}{Theorem}
\newtheorem{theorem}{Theorem}

\newtheorem*{corollary*}{Corollary}

\usepackage{alltt}  %

 \newcolumntype{Z}{ >{\centering\arraybackslash}X }

\usepackage[autolanguage, np]{numprint}

\begin{document}

%
%

%
%
\textbf{A Stochastic Record-Value Approach to Global Simulation Optimization}

\begin{center}
    \textit{Rohan Rele}
\\

\textit{
Committee:
Zelda Zabinsky , 
Aleksandr Aravkin , 
Giulia Pedrielli }

\textit{Master of Science:}
A thesis submitted in partial fulfillment of the requirements for the degree of Applied Mathematics, 2021 
\end{center}

%
%

%
%

\setcounter{page}{-1}
\abstract{\indent Black-box optimization is ubiquitous in machine learning, operations research and engineering simulation. Black-box optimization algorithms typically do not assume structural information about the objective function and thus must make use of stochastic information to achieve statistical convergence to a globally optimal solution. \\
\indent One such class of methods is multi-start algorithms which use a probabilistic criteria to: determine when to stop a single run of an iterative optimization algorithm, also called an \textit{inner search}, when to perform a restart, or \textit{outer search}, and when to terminate the entire algorithm. Zabinsky, Bulger \& Khompatraporn introduced a record-value theoretic multi-start framework called Dynamic Multi-start Sequential Search (DMSS). We observe that DMSS performs poorly when the inner search method is a deterministic gradient-based search. \\
\indent In this thesis, we present an algorithmic modification to DMSS and empirically show that the Revised DMSS (RDMSS) algorithm can outperform DMSS in gradient-based settings for a broad class of objective test functions. We give a theoretical analysis of a stochastic process that was constructed specifically as an inner search stopping criteria within RDMSS. We discuss computational considerations of the RDMSS algorithm. Finally, we present numerical results to determine its effectiveness. 
 }
 
%
%
\tableofcontents
\listoffigures
 \comment{
%
%

\section*{Glossary}      
\addcontentsline{toc}{chapter}{Glossary}
\thispagestyle{plain}
\begin{glossary}

\end{glossary}
 }
%
%

%
%
%

%
%


\section{Introduction}

Mathematical optimization traditionally seeks to minimize or maximize an objective function subject to some constraints. In unconstrained optimization, properties like smoothness or convexity of the objective function facilitate the derivation of optimality conditions that can be used to guarantee analytic convergence to locally or globally optimal solutions. 

In black-box optimization, such properties are not generally available nor are they practical to determine globally. Algorithms like Simulated Annealing \cite{annealing} converge in probability to a global optimum for continuous objectives because they leverage distributional properties of Markov transition probabilities or \textit{restart} probabilities. The usage of probabilistic restart conditions persists throughout this thesis and motivates the primary class of global, black-box optimization algorithms that we consider, known as \textit{multi-start algorithms}. 

The goal of global optimization is to obtain a \textit{global solution} to an objective function $f$. A general optimization problem is given by: 
\begin{equation}
    \min_{x \in \mathbb{X}} f(x)
\end{equation} where we define a global solution (sometimes referred to as an \textit{optimum}) as a value $x^* \in \mathbb{X}$ such that $f(x^*) \leq f(x)$ for all $x \in \mathbb{X}$. Another related definition is that of a \textit{local optimum} or \textit{solution}. A local optimum is defined as an $\Tilde{x} \in \mathbb{X}$ such that $f(\Tilde{x}) < f(x)$ for all $x \in \mathcal{N}(\Tilde{x}) \subset \mathbb{X}$ where $\mathcal{N}(\Tilde{x})$ denotes a neighborhood around $\Tilde{x}.$ 

The existence of global optima under specific conditions is the subject of many texts in the mathematical optimization community \cite{Horst2000IntroductionTG}. Here we state some general theorems in order to proceed with the fewest possible assumptions about Problem (1.1).  

\begin{theorem}[Weierstrass]
If the objective function $f$ is continuous and $\mathbb{X}$ is closed and bounded, then $x^* \in \mathbb{X}$ exists.  
\end{theorem}

This is a fundamental theorem in global optimization \cite{Horst1995HandbookOG}. It implies that if the feasible region of a continuous function $f$ is real-valued, its compactness is sufficient to achieve a globally optimal solution. A corollary of this can be stated for a \textit{closed} feasible region with an additional condition on $f$. 

\begin{corollary*}
If the objective function $f$ is continuous and \textit{coercive} i.e., $$ \lim_{||x|| \rightarrow \infty} f(x) = +\infty $$ on a closed $\mathbb{X}$, then $x^* \in \mathbb{X}$ exists. 
\end{corollary*}

What is notable about these two results is their generality, which allows for application to black-box functions. Many real systems that are modeled by a computational oracle can only be calculated on closed or compact domains. Coercivity is also a realistic property in that simulated functions generally will not admit upper/lower bounds. With these results in mind, we make the following key assumption that allows many algorithms to achieve asymptotic convergence: 

\begin{quote}
    \textit{The black-box function $f$ is equipped with an oracle $f(.)$ such that we can calculate a value $f(x)$ for all $x \in \mathbb{X}.$}
\end{quote}

\comment{We aim to adaptively search for the global optimum $\textbf{x}^*$ via a nested combination of deterministic and/or stochastic algorithms. We synchronize the execution and interaction of these algorithms in a specific manner to achieve statistical performance metrics. The performance metrics for the combination of algorithms include asymptotic convergence to $f^*$, and statistical confidence bounds on achieving $\epsilon-$optimality within a fixed number of oracle calls. }

Given the above, we state the formal optimization setting that we remain in henceforth:

\begin{equation}
    \min_{x \in \mathcal{X}} f(x)
\end{equation}
where $f: \mathcal{X} \subset \mathbb{R}^d \rightarrow \mathbb{R}$ has a compact domain $\mathcal{X}$.  We refer to the global minimum of the function $f$ with respect to the domain $\mathcal{X}$ as $x^* \in \displaystyle \arg\min_{x \in \mathcal{X}} f(x)$ with corresponding minimum value $f^* = f(x^*)$. \\

\comment{The function $f$ is black box which means that we are assumed to be given a computationally expensive black-box oracle for $f(x)$ for all $x \in \mathcal{X}$. Note that since $f$ is black-box, i.e. it may be non-linear, non-convex or non-smooth, we do not assume any additional structural properties.}

The outline of this thesis is as follows. In Chapter 2 we discuss adaptive random search and multi-start algorithms. We focus on Dynamic Multi-start Sequential Search (DMSS) as an algorithmic framework but also discuss its underlying conceptual algorithm Hesitant Adaptive Search with Power-Law Improvement Distribution (HASPLID) \cite{HASPLID}. In Chapter 3 we present the Revised DMSS (RDMSS) algorithm along with an analysis of the stochastic record-improvement slope process, which is built off of the assumptions of HASPLID. Chapter 4 presents the numerical results and experimental evidence to examine the choice of RDMSS over DMSS for gradient-based inner search settings. Lastly, we discuss future directions in Chapter 5.

 
\section{Background}

\begin{figure}
    \centering
    \includegraphics[width=0.5\textwidth]{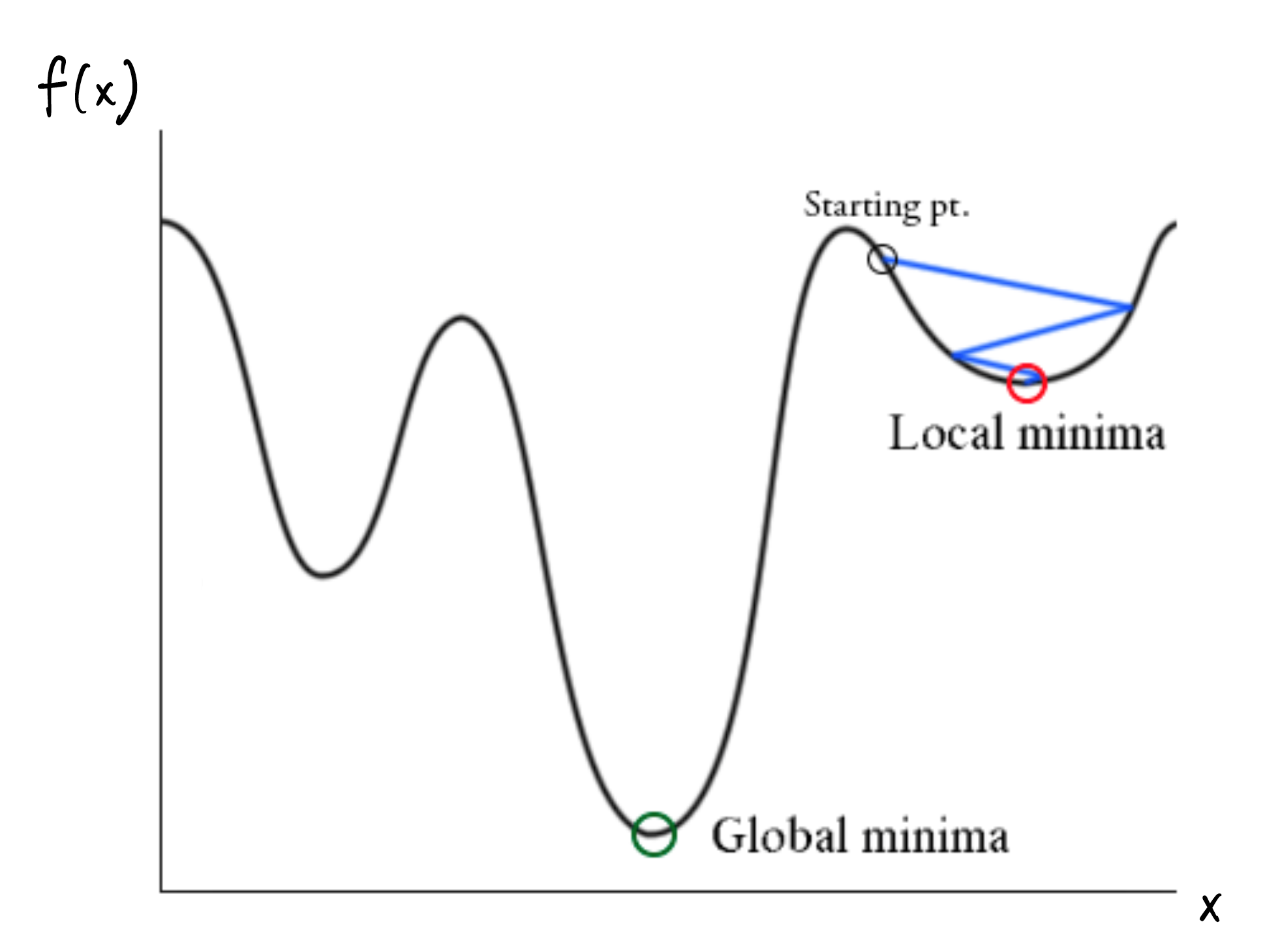}
    \caption{Classical gradient descent algorithm getting ``trapped" at a local optimum \cite{grad_trapped}.}
    \label{fig:grad_descent_trapped}
\end{figure}

Unconstrained nonlinear programming approaches such as gradient descent, conjugate gradient methods and trust region methods are examples of \textit{local} optimization techniques or \textit{local searches} \cite{local_opt_review}. What distinguishes them from \textit{global} optimization techniques is that they are not guaranteed to achieve global optimality under the most general conditions. \\
\indent Gradient descent, one of the most commonly-used first order optimization methods in scientific computing, is said to have achieved an optimal solution if, for some smooth $f$, $$\nabla f(\Tilde{x}) = 0 $$ for some $\Tilde{x}$ that lies in a general domain $\mathbb{X}$. However, without additional information such as convexity of $f$, there is not sufficient information to determine if $\Tilde{x} = x^*$, namely, if $\Tilde{x}$ is a globally optimal solution. Many times in practice, it is not and it is for this reason that gradient-based methods suffer from the phenomena of being ``trapped" in local optima, hence the name local search. See Figure \ref{fig:grad_descent_trapped}. \\

\section{Overview of Multi-Start}

\indent To remedy this, the technique of \textit{randomly restarting} a local search has been studied extensively \cite{glob_opt_based_local_search}. Upon completion of a local search, the modified optimization algorithm equipped with restart functionality employs a stochastic criterion to initiate a restart, at which point the local search will resume again from a potentially unexplored initial point in the domain. \\
\indent A \textit{multi-start} framework \cite{Marti2003} formalizes the idea of performing successive local optimization algorithms with different restart locations as initial starting points. There are four key components that characterize multi-starts. These are the: \\

\begin{enumerate}
    \item[(i)] outer loop
    \item[(ii)] outer search 
    \item[(iii)] inner loop
    \item[(iv)] and inner search.
\end{enumerate}

\begin{figure}
     \centering
     \includegraphics[width=0.8\textwidth]{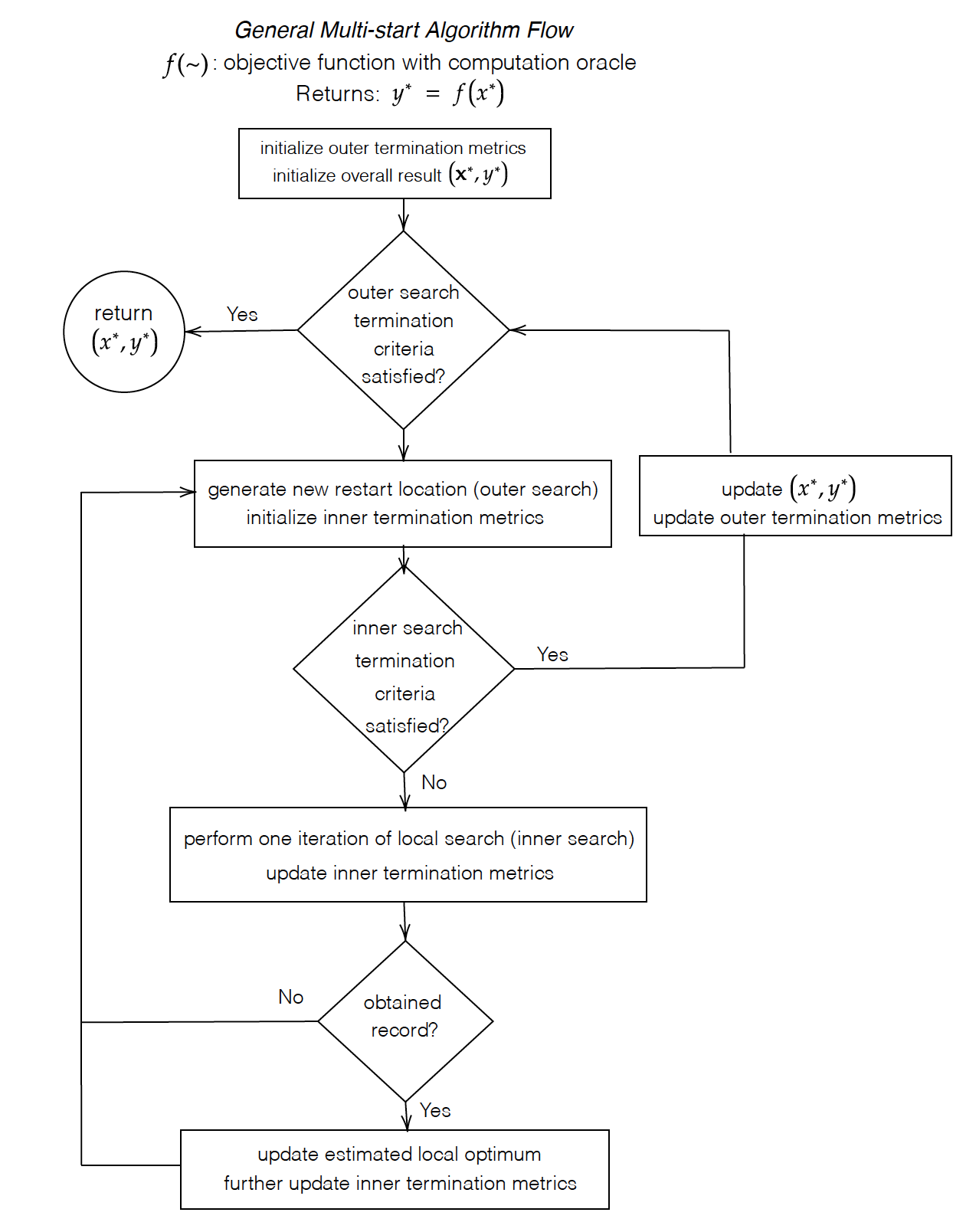}
     \caption{General Multi-Start Algorithm}
     \label{fig:general_flow}
 \end{figure}

\indent The outer loop serves as a global termination criterion verifier. In other words, it is the particular statement within the algorithm that decides whether or not it has found a sufficiently close estimate of the global optimum. It is responsible for iteratively performing outer searches. The outer search is a sampling method that stochastically regulates \textit{where} in the domain the algorithm should search next. In other words, the outer search generates the initial point for each \textit{independent} run of the particular local optimization method that the algorithm employs. Most multi-start algorithms assume independence of local searches, however some methods explore dependence, for example \cite{SOAR}. The importance of the independence between outer searches cannot be over-stated, it is a key assumption in this thesis which we shall return to later. \\
\indent The inner search and inner loop are components of the local optimization method. In the example of gradient descent as the local search method, consider a single update of the sequence of estimated optimal locations $(x_n)_{n=1}^N$ of a smooth $f$ with step-size $\gamma$: 
$$x_{n+1} = x_n - \gamma\nabla f(x_n). $$
A single such update is what is referred to as an \textit{inner search} while the sequence of inner searches is terminated subject to a condition imposed by the \textit{inner loop}, analogously to the relationship of the outer search and its loop. More generally, an inner search is simply one search iteration of a local optimization method given some initial point and the inner loop determines when a local termination criterion has been reached. We tie the above concepts together with a high-level flowchart (Figure \ref{fig:general_flow}) and associated pseudo-code (Algorithm \ref{alg_1}). \\
\indent Returning to the example of gradient descent once more, we can see that $f(x_{n+1}) \leq f(x_n)$ by construction. The optimality condition makes the inequality tight and, as such, no improvements can be made to the estimated optimum. Another way to frame this is in the language of \textit{records}. A record value is simply an improvement over the running best estimate of an optimal point. This may refer to an estimated local optimum or an estimated global optimum. In the case of gradient descent, every $x_n$ is a record since $f(x_n) < f(x_{n-1})$ for all $n < N$ until the optimality condition is achieved, at which point $x_N = \Tilde{x}$. The same property holds for all gradient-based local searches.\\ 
\indent The analysis of records is more interesting for global methods. Successive local searches performed by the outer loop may or may not yield continuously improving results. It is for this reason that the outer search, the stochastic sampling method, is often used to derive a probabilistic global termination criteria. The independent samples allow us to perform a statistical analysis of ``when" records should appear. This theory, called record-value theory, is intimately related to order statistics \cite{Records},\cite{gupta_1984},\cite{nagaraja_nevzorov_1997}. \\
 \begin{algorithm}
 \KwData{$f:$ objective function with computational oracle}
 \KwResult{$(\textbf{x}^*, \textbf{y}^*)$: global minimum $\textbf{y}^* = f(\textbf{x}^*)$ and associated location $\textbf{x}^*$}
 \noindent\rule{12cm}{0.4pt}\\ 
\textbf{Algorithm}($f$): \\
initalize outer termination metrics \\
initialize $(\textbf{x}^*, \textbf{y}^*)$ \\
\While{outer termination criteria not satisfied}{
generate new restart location $\leftarrow$ perform outer search \\
initialize inner termination metrics \\
\While{inner termination criteria not satisfied}
{
obtain search iterate $\leftarrow$ perform inner search \\
update inner termination metrics \\
\If{a record was obtained}
    {update estimated local optimum \\
    further update inner termination metrics \\
    }
}
update $(\textbf{x}^*, \textbf{y}^*)$ \\
update outer termination metrics
}
\Return $(\textbf{x}^*, \textbf{y}^*)$ \\

 \caption{General Multi-Start Algorithm}
 \label{alg_1}
 \end{algorithm}
 \indent Before discussing the main example of a multi-start framework which we will be dissecting in the remainder of the thesis, we present common considerations when developing or analyzing multi-starts. That is, how to balance \textit{exploration} vs. \textit{exploitation}. Exploration refers to efficient discovery of the sample space. In black-box settings, even with compactness or closure of the domain, the issue of discovering new sample points is routinely emphasized in developing computational solvers. To add further complexity, we want to consider any potential underlying structure and \textit{exploit} the information that is being gathered during each local search. In other words, we would like our global optimization algorithm to cover as much ``ground" on the domain as possible, since a global optimum may lie anywhere on it, while simultaneously making intelligent use of past discovery. \\ 
 \indent An illustrative example of an intelligent restart strategy is the Stochastic Optimization for Adaptive Restart (SOAR) framework \cite{SOAR}. SOAR \textit{exploits} the previously collected information by making use of a Gaussian Process as a surrogate model for its outer search method. Thus, the restart location is chosen via a posterior distribution that is updated upon every function evaluation throughout the algorithm. This yields a statistical basis for how the algorithm chooses to progress. In order to also use the surrogate model as a means of \textit{exploration}, a criterion is built into the sampler that fixes the algorithm's subsequent search space to the set of un-sampled points. Meaning, if the algorithm determines that the estimated Gaussian Process on any particular iteration is not expected to yield a ``successful" local search at the distributionally-generated point, then the algorithm will restart at a completely new, unexplored location and resume Gaussian Process updates from there. Here, we use the word ``successful" to mean whether or not the local search yields a record and the word ``expected" in both its colloquial and probabilistic sense. \\
\indent Some analytic considerations of multi-start include (1) the conditions upon which to perform a restart, (2) the global termination criteria, i.e., when to \textit{stop} performing restarts, (3) the number of iterations taken by the local search per restart, (5) the number of iterations required to first achieve a sufficiently close estimate of the global optimum and (6) the total number of restarts. Unlike SOAR, which takes a fixed computational budget and decrements it until it has been fully expended, Dynamic Multi-start Sequential Search (DMSS) \cite{HASPLID} addresses (1), (2) and (3), choosing to use an adaptive computational expense metric to progress successive local searches. The difference is that DMSS assumes independent restart sample points and cannot make use of distributionally learned information to sample in the way that SOAR's Gaussian Process did. We revisit these considerations closely in the following section.

\section{Dynamic Multi-start Sequential Search}

Dynamic Multi-start Sequential Search (DMSS) is a specific multi-start framework that exploits results from record-value theory in order to define stochastic metrics that limit the computational expense of finding the global optimum. Before performing an in-depth analysis of these metrics, we note that DMSS is modeled by Hesitant Adaptive Search with Power-Law Improvement Distribution (HASPLID) \cite{HASPLID}. We summarize the key results of the conceptual algorithm HASPLID in order to motivate the subsequent analysis of DMSS and its variant, the Revised DMSS algorithm. 

\subsection{Hesitant Adaptive Search with Power-Law Improvement Distribution}

\medskip

The stopping and restarting conditions for DMSS are built on a parametrized conceptual model called the Hesitant Adaptive Search with Power-Law Improvment Distribution (HASPLID) algorithm. HASPLID formulates criteria for stopping a single run of an executed search for an optimum and determines whether to restart another run or terminate the whole algorithm. A key quantity in HASPLID is the range distribution $\rho$, which is an implicitly defined measure that is written in terms of the sampling distribution $\mu$ on $\mathcal{X}$: 
$$\rho(T) = \mu(f^{-1}(T)) $$ for $T \in \mathcal{B}(\mathbb{R})$ and its CDF $p(y) = \mu(f^{-1}((-\infty, y]))).$\\
\indent We assume the distribution $\rho$ to be continuous and introduce two parameters: $\alpha \in [0,1]$ which controls the difficulty of finding improvements and $\lambda \in \mathbb{R}^+$, which controls the distribution of the quality of improvements found. As in \cite{HASPLID}, we consider a power-law transformation $\rho^{(\lambda)}$ of the range distribution $\rho$ for an arbitrary termination set $T \subseteq \mathbb{R}$: 
$$\rho^{(\lambda)}(T) = \int_{z\in T} d(p(z))^\lambda = \lambda \int_{z\in T} (p(z))^{\lambda - 1} dp(z). $$ The relationships between CDF $p^{(\lambda)}$ of $\rho^{(\lambda)}$ and CDF $p$ of $\rho$ still hold as defined above but are now specifically given by: $$p^{(\lambda)}(y) = (p(y))^\lambda $$ The normalized restriction of $\rho^{(\lambda)}$ to the left half-line $(-\infty, y] $ is given by: $$\rho_y^{(\lambda)}(T) = \lambda(p(y))^{-\lambda} \int_{z \in T \bigcap (-\infty,y]} (p(z))^{\lambda-1} dp(z) $$ with CDF $$p_{y_1}^{(\lambda)}(y_2) = \left( \frac{p(y_2)}{p(y_1)} \right)^\lambda $$
for $y_1 \leq y_2 < y$ where $y_1, y_2 \in \mathcal{X}$. 

\medskip 

\noindent The parametrized algorithm HASPLID($\alpha, \lambda$; $\rho$) is given by the following pseudocode: \\

\noindent \textbf{HASPLID}($\alpha, \lambda$; $\rho$) \textit{c.f.} \cite{HASPLID}
\begin{quote}
    \begin{enumerate}
        \item[\textbf{Step 0.}] Set $j = 0$. Sample $Y_0$ according to $\rho^{(\lambda)}$. 
        \item[\textbf{Step 1.}] With probability $(p(Y_j))^\alpha$, sample $Y_{j+1}$ according to $\rho_{Y_j}^{(\lambda)}$. With probability $1 - (p(Y_j))^\alpha$, set $Y_{j+1} = Y_j$. Note that this construction yields: $Y_0 \geq Y_1 \geq \cdots \geq Y_j$ in the case of minimization. 
        \item[\textbf{Step 2.}] If a stopping criteria is met, stop; otherwise, increment $j$ and return to \textbf{Step 1}. 
    \end{enumerate}
\end{quote}

\begin{figure}[!h]
    \centering
    \includegraphics[width=\textwidth]{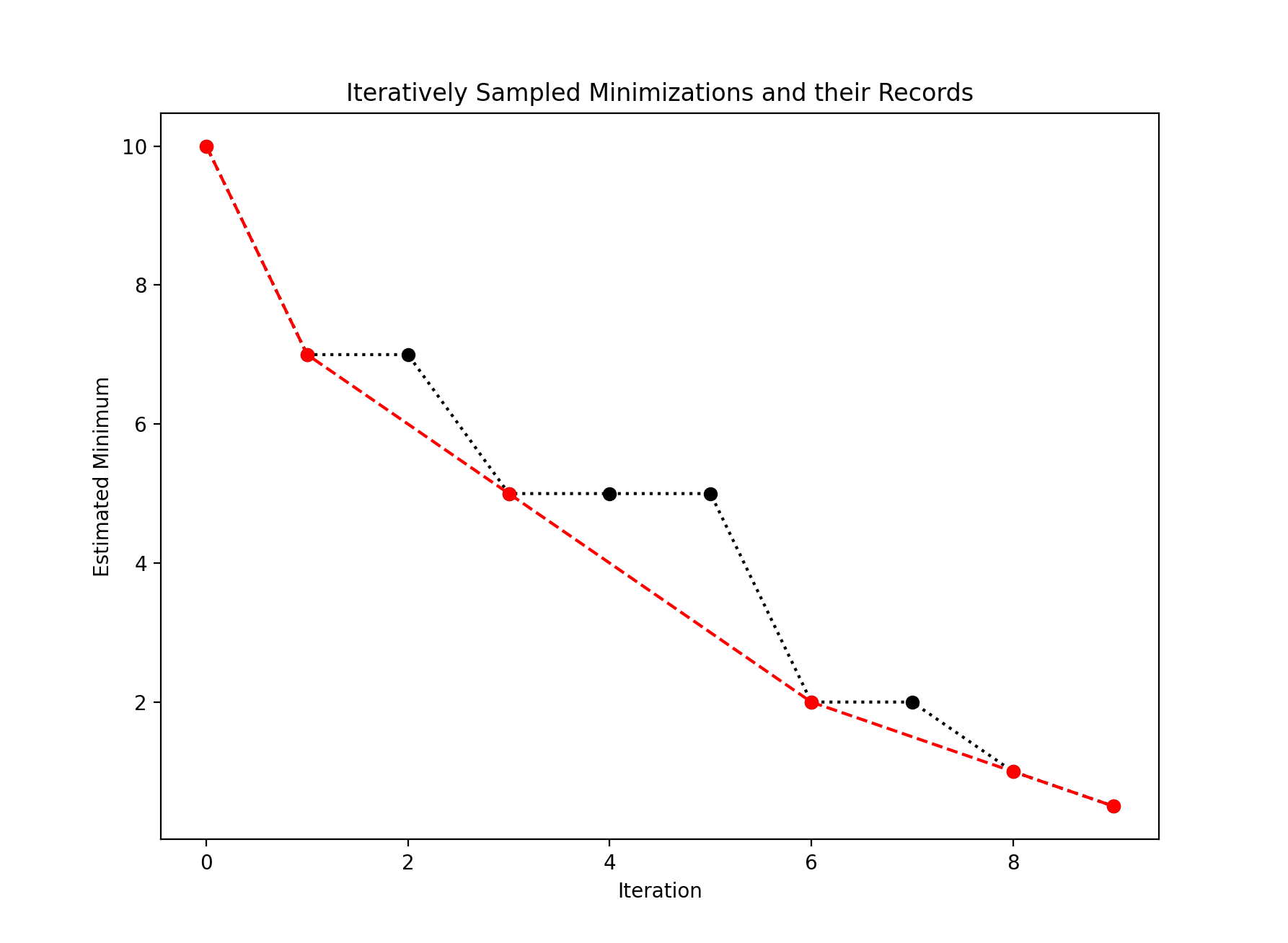}
    \caption{Graphical representation of inner search iterates and their record values. Every point is the result of an inner search but only the red points are records. For example, $Y_4 = 5$ but $Y_{R(5)} = 1$ and $R(5) = 8$.}
    \label{fig:Record_Example}
\end{figure}

\indent The weak monotonicity of the random variables $(Y_j)$ is significant. HASPLID is a general adaptive search algorithm that decides, at each sampling iteration, whether to sample in the set of improving points with an associated probability $(p(y))^\alpha$, called the \textit{bettering probability} \cite{GlobalOpt}, or remain at its currently sampled iterate with probability $1-(p(y))^\alpha$. This implies that at each iteration $j$ of a HAS algorithm, $Y_{j+1}$ will either be a record $Y_{j+1} < Y_j$ with probability $(p(y))^\alpha$ or $Y_{j+1} = Y_j$ with probability $1-(p(y))^\alpha$. \\
\indent We state this formally and introduce notation to distinguish between records and raw search iterates: $(Y_j)_{j\in \mathbb{N}}$ forms a non-increasing Markov Chain in which no result of any given search is greater than the running minimum over the set of all random searches. By contrast, the collection of \textit{records} $(Y_{R(k)})_{k,R(k) \in \mathbb{N}}$ is a strictly monotone decreasing Markov Chain that represents every improving result over the set of all terminated random searches. In this notation, $Y_j$ is a record denoted $Y_{R(k)} := Y_j$ with $R(k) = j$ if $Y_{R(k-1)} > Y_j$ for $j > R(k-1)$. See Figure \ref{fig:Record_Example} for a graphical depiction and example. 

\indent We now characterize some of the key behaviors of HASPLID which underlie DMSS: \\ 

\comment{
\noindent \textbf{Proposition 1} (\textit{c.f.} see \cite{HASPLID}) The iterates $(Y_j)$ of HASPLID($\alpha,\lambda$;$\rho$) are stochastically equivalent to the iterates $(\Tilde{Y}_j)$ of HASPLID($\alpha/\lambda,1$;$\rho$). \\
}
\noindent \textbf{Proposition 2} (\textit{c.f.} see \cite{HASPLID}) Let $N(y)$ denote the number of records obtained by HASPLID($\alpha,\lambda$;$\rho$) before obtaining a value of $y$ or better. Then $N(y) \sim \text{Poi}(-\lambda \log p(y)).$ \\

\noindent \textbf{Proposition 3} (\textit{c.f.} see \cite{HASPLID}) The probability $$\mathbb{P}(Y_{R(k)} > y) = 1 - (p(y))^\lambda \sum_{s=0}^{k-1} \frac{(-\lambda \log p(y))^s}{s!} = G(k, -\lambda \log p(y)). $$ where $G(n,x) = 1 - e^{-x}\sum_{s=0}^{n-1} x^s/s!$ denotes the incomplete gamma function. 

\noindent \textbf{Proposition 5} (\textit{c.f.} see \cite{HASPLID}) The probability of obtaining $k$ records in the first $j$ iterates of HASPLID$(\alpha, \lambda; \rho)$ is: $$(\lambda/\alpha)^{j-1}|s(j,k)| \cdot \frac{\Gamma(1 + \lambda/\alpha)}{\Gamma(j + \lambda/\alpha)} $$ where $s(j,k)$ is a Stirling number of the first kind and $\Gamma(x)$ denotes the \textit{complete} gamma function. 

\noindent \textbf{Proposition 6} (\textit{c.f.} see \cite{HASPLID}) The expected number of records found by HASPLID($\alpha, \lambda; \rho$) in the first $j$ iterates is: $$\frac{\lambda}{\alpha}\left( \psi\left(j + \frac{\lambda}{\alpha}\right) - \psi\left(\frac{\lambda}{\alpha} \right) \right) $$ where $\psi$ is the digamma function. 

\noindent \textbf{Proposition 7} (\textit{c.f.} see \cite{HASPLID}) The probability that $R$ independent runs of HASPLID$(\alpha, \lambda; \rho)$, with $k_1, k_2, \cdots, k_R$ records obtained respectively, never approach the $\epsilon-$target region around the global optimum is:

\begin{equation}
    \prod_{r=1}^R G(k_r, -\lambda \log \epsilon). 
\end{equation}
 
\subsection{Outer Loop (DMSS)}
We begin our discussion of DMSS from the ``top-down". This section is dedicated to the characterization of DMSS' outer termination metrics and global termination criterion. Like all multi-start algorithms, the outer loop addresses the analytic consideration of ``when to stop performing restarts", this is implemented by the global termination condition. DMSS introduces a second consideration, the ``total number of restarts", into the global termination condition. In particular, the DMSS global termination condition leverages past restarts in deciding whether or not the globally estimated optimum is sufficient. \\
\indent Before defining the condition, we consider the information available to the algorithm upon start. The conceptual HASPLID algorithm takes two parameters, $\alpha \in [0,1]$ and $\lambda \in \mathbb{R}^+$ from which DMSS inherits only the parameter $\alpha$ and introduces two new parameters $0 < \epsilon \ll 1$ and $0 < \delta \ll 1$. The parameter $\epsilon$ specifies the size of the $\epsilon-$target region, defined as the $(2\epsilon)^d$-sized hypercube around the global optimum where $d = \dim(\mathcal{X})$ is the dimension of the domain of the objective function $f: \mathcal{X} \rightarrow \mathbb{R}$. We will return to $\delta$ shortly. \\
\indent DMSS invokes Proposition 7 in defining its primary outer termination metric:

\begin{equation*}
    p_{\text{FAIL}} = \prod_{r=1}^R G(k_r, -\lambda \log \epsilon) 
\end{equation*}

\noindent where $G$ is the incomplete gamma function. The $\delta$ parameter simply serves as a user-defined upper bound for the $p_{\text{FAIL}}$, the probability that DMSS fails to achieve a sufficiently close estimate of the global optimum. Thus the global termination criterion is nothing more than $p_{\text{FAIL}} < \delta$. The form of this probability in (2.1) makes it clear why \textit{independence} was needed as an assumption when motivating multi-starts; without it, the global termination criteria could not be defined this way. 

\indent Note the absence of $p(y)$ in the $p_{\text{FAIL}}$ metric. Proposition 5 makes use of a limiting property of Stirling numbers and characterization of HASPLID iterates as exponentially-distributed record statistics to vanish the distribution $p(y)$ that appears in Propositions 1-3. 

\indent One computational issue is that DMSS cannot calculate explicitly $p_{\text{FAIL}}$ without $\lambda$. To combat this \cite{HASPLID} invokes Proposition 1 of \cite{HASPLID} and works instead with the ratio $\lambda/\alpha$. They define $\zeta := \frac{\lambda}{\alpha}$ and solve the following maximum likelihood equation for $\zeta$ on each run 
\begin{equation}
    \sum_{r=1}^R(k_r - 1) + \zeta\left(R\psi(1 + \zeta) - \sum_{r=1}^R \psi(j_r + \zeta) \right) = 0
\end{equation}
which is derived from Proposition 5 of \cite{HASPLID}. This allows us to replace $\lambda$ with $\alpha\zeta$ in (2.1).

\subsection{Outer Search (DMSS)}

The power of DMSS lies in its outer/inner search flexibility. The outer search in particular is responsible for generating restart locations which are fed to the first inner search within a particular run. These sampling locations only require independence from one another. In theory, the sampling locations need not be identically distributed; however, in practice, i.i.d. uniform random samples are the most straightforward to implement. In all analyses going forward, we fix the outer search distribution to be i.i.d uniform random. \\

\subsection{Inner Loop (DMSS)}

The inner loop of DMSS, like its outer loop, couples two analytic considerations to formulate the termination criterion. Like all multi-start algorithms, the inner loop addresses when to perform a restart but it does so using the \textit{number of records} found on a particular run. We note that Proposition 2 \cite{HASPLID} states that the number of records required to obtain a particular value $y$ is distributed as a Poisson random variable with parameter $-\lambda \log p(y)$. It can be shown \cite{nagaraja_nevzorov_1997} that this implies the ``time" between records or \textit{inter-record time} is distributed as a \textit{geometric} random variable with parameter $(p(y))^\alpha$. In other words, $$R(k+1) - R(k) \sim \text{Geo}((p(y)^\alpha) $$ where $Y_{R(k)} = y$. Here we use \textit{time} to mean ``number of raw search iterates". This fact, along with independence between record values and inter-record times, allows the use of Propositions 5 and 6 to define the primary inner termination metric 

\begin{equation}
    n_{\text{RECORD}} := \frac{\lambda}{\alpha}\left( \psi\left(j + \frac{\lambda}{\alpha}\right) - \psi\left(\frac{\lambda}{\alpha} \right) \right) 
\end{equation}

\noindent which represents the expected number of raw search iterates prior to achieving a record. Note that this value depends on $j$, the number of raw search iterates, so this value updates as more and more inner searches are performed. Intuitively, we would like stop performing inner searches when the number of raw search iterates exceeds $n_{\text{RECORD}}$ meaning that it is taking longer than expected to find a record. This is exactly the inner (local) termination criterion: $j_R < n_{\text{RECORD}}$ where $j_R$ is the number of raw search iterates on run $R$. \\

\indent Like the outer loop, we note several computational issues. DMSS continues making use of $\alpha\zeta$ in place of $\lambda$. Additionally, we must now devise a way to distinguish between the number of records in run $R$ and the number of raw search iterates in run $R$. We use the indexed variables $k_R$ and $j_R$ respectively. Clearly $j_R \geq k_R$.

\subsection{Inner Search (DMSS)}

The inner search of DMSS is once again flexible. In \cite{HASPLID}, the numerical results presented fix an elitist random ball walk as the inner search method, which is stochastic. However, this can easily be modified to be a deterministic search such as gradient descent. In the following chapter, we discuss a particular consideration of a deterministic, gradient-based inner search and how it motivates the Revised DMSS algorithm. 

\begin{figure}
     \centering
     \includegraphics[width=0.7\textwidth]{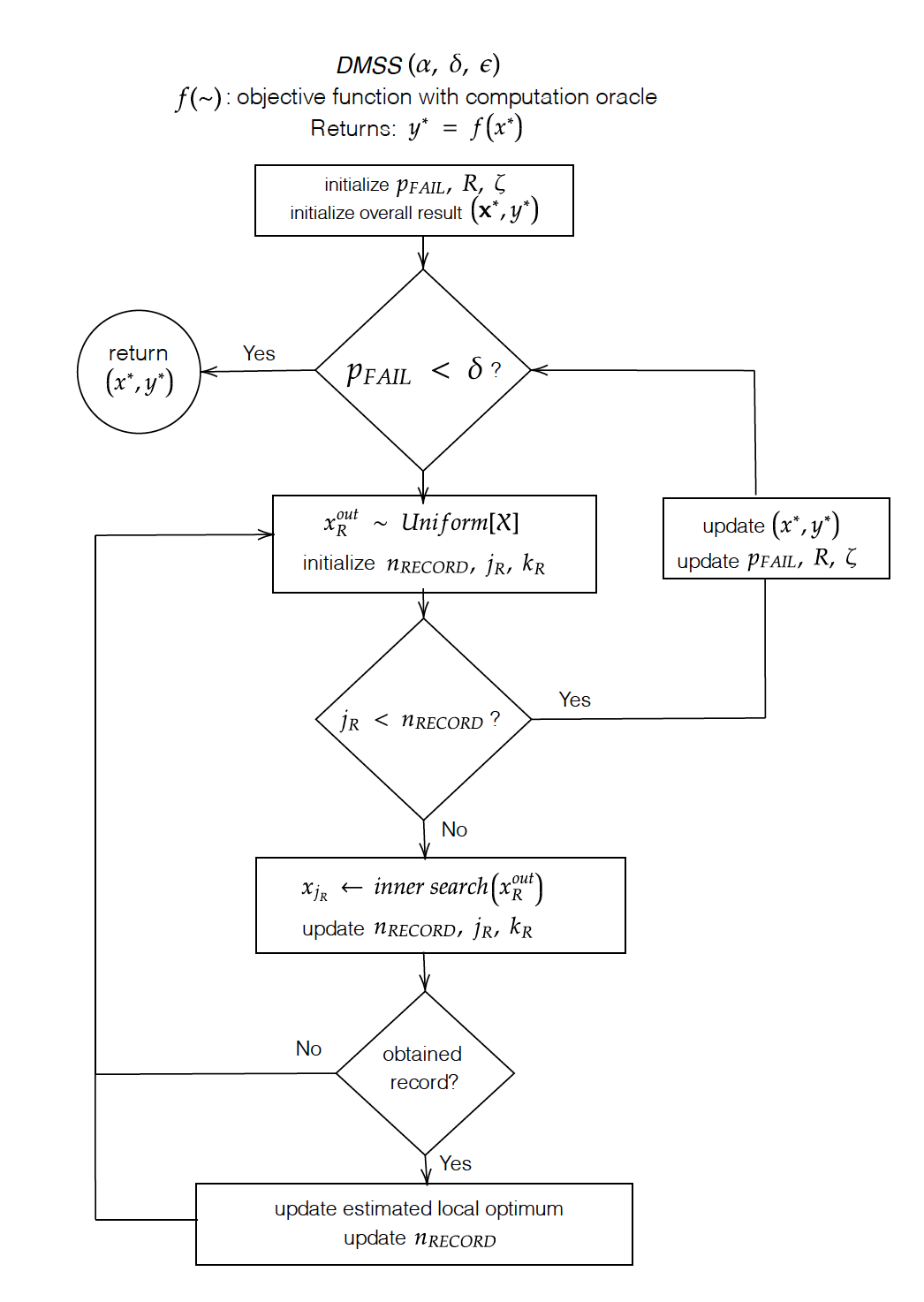}
     \caption{DMSS($\alpha, \delta, \epsilon$)}
     \label{fig:dmss_flow}
 \end{figure}

\indent Putting these components together yields the full DMSS algorithm. See Figure \ref{fig:dmss_flow} and Algorithm 2. 

 \begin{algorithm}
 \KwData{$f:$ objective function with computational oracle}
 \KwResult{$(\textbf{x}^*, \textbf{y}^*)$: global minimum $\textbf{y}^* = f(\textbf{x}^*)$ and associated location $\textbf{x}^*$}
 \noindent\rule{12cm}{0.4pt}\\ 
\textbf{Algorithm}($f$): \\
initalize $p_{\text{FAIL}}, R, \zeta$ \\
initialize $(\textbf{x}^*, \textbf{y}^*)$ \\
\While{$p_{\text{FAIL}} < \delta$}{
$y_R^{\text{out}} \leftarrow Uniform[\mathcal{X}]$ \\
initialize $n_{\text{RECORD}}, j_R, k_R$ \\
\While{$j_R < n_{\text{RECORD}}$}
{
$x_{j_R}$ $\leftarrow$ inner search$(x_R^{\text{out}})$ \\
update $n_{\text{RECORD}}, j_R, k_R$ \\
\If{a record was obtained}
    {update estimated local optimum \\
    update $n_{\text{RECORD}}, k_R$ \\
    }
}
update $(\textbf{x}^*, \textbf{y}^*)$ \\
update $p_{\text{FAIL}}, R, \zeta$
}
\Return $(\textbf{x}^*, \textbf{y}^*)$ \\
\caption{DMSS}
\end{algorithm}


\section{Revised Dynamic Multi-start Sequential Search}

Preliminary observations made during experimental tests showed that DMSS did not perform well when the inner search was a determinstic, gradient-based search. The reason for this is the inner termination criteria. Recall that a deterministic, gradient-based optimization method obtains a record \textit{on every iteration except its last} by definition. This means that $k_R$ and $j_R$ of DMSS are being updated simultaneously on every iteration of the inner loop until a local optimum has been found. From an algorithmic standpoint, the local search seems to be doing extremely well and thus, by the construction of $n_{\text{RECORD}}$, the inner termination metric of DMSS, it remains relatively large. But when the local optimum is found, $k_R$ will stop updating and $j_R$ will need to iterate many more times before exceeding $n_{\text{RECORD}}$. This results in a waste of computational expenditure and mimics the ``trapping" phenomena that make gradient-based searches unsuitable for global optimization of multimodal test functions. \\
\indent The Revised DMSS (RDMSS) algorithm was designed primarily to address the issue discussed above. Our algorithm was formulated in a manner that leverages HASPLID and continues to take advantage of the theoretical results that characterize number of record values per run. We accomplish this by introducing a \textit{second} inner termination metric that works in tandem with the original metric, $n_{\text{RECORD}}$, to sharpen the criteria upon which a run is terminated. While $n_{\text{RECORD}}$ regulates the number of raw search iterates, the new metric regulates the \textit{difference in magnitude between record values}. \\
\indent To provide experimental control and ultimately determine whether or not RDMSS outperforms DMSS in gradient-based settings, we fix the inner searches of both algorithms to be an arbitrary determinstic, gradient-based search henceforth. For the experiments performed in Chapter 4, this method is chosen to be the Newton Conjugate-Gradient method.  \\

\section{Outer Loop (RDMSS)}
The outer loop of Revised DMSS continues using the outer termination metrics/criteria of DMSS. This was done to leverage the HASPLID-specific results that were derived in \cite{HASPLID}. The re-use of the independence assumption allows us to use $p_{\text{FAIL}}$ for the outer termination criteria without modification. Hence, the outer loops of RDMSS and DMSS are identical. 

\section{Outer Search (RDMSS)}
The outer search of Revised DMSS is also unchanged from its predecessor. The reason for this is, once again, the independence assumption that results in the $p_{\text{FAIL}}$ metric. Relaxing this assumption means re-examining the question of \textit{exploration} vs. \textit{exploitation} and possibly losing the ability to use the expected number of records to formulate our inner and outer termination criteria. Thus, for RDMSS, we continue to assume i.i.d uniform random restart locations across the domain.

\section{Inner Loop (RDMSS)}

The inner search termination criterion is where RDMSS deviates from DMSS. Our new criterion is built upon the following definition: \\
    \noindent \textbf{Def. } The \textit{slope process} is a collection of random variables $(S_k)$ with a common distribution. Each $S_k$ is given by the ratio 
    \begin{equation}
        \displaystyle S_k := \frac{Y_{R(k)} - Y_{R(k+1)}}{R(k+1) - R(k)}.
    \end{equation}

This definition relies heavily on the assumptions and implications of HASPLID. We are provided an extensive analysis of the \textit{time between records} or the number of sampling iterations between record $k$ and record $k+1$ for all $k \geq 1$ in \cite{HASPLID}. Further, \cite{Records} explains that for the \textit{classical} record model, the raw sampling iterates $(Y_j)$ are i.i.d exponential random variables with parameter $1$. Generalizing slightly, we let $(Y_j) \sim \text{exp}(\Lambda).$ Due to the memory-less property of the exponential distribution, the sequence $(Y_{R(k)} - Y_{R(k+1)})_{k \geq 0} \sim G(k+1,\Lambda)$ for all $k$ where $G$ denotes the incomplete gamma function. Note that under the assumption that $(Y_j)$ has a continuous common CDF (in the case of HASPLID, this is defined as $p(y)$), the survival function for the $(k+1)$-st record value is given by $$ \mathbb{P}(Y_{R(k+1)} > y^*) = 1 - e^{-y^*} \sum_{s = 0}^{k} (y^*)^s/s! \text{ for } y^* > 0.$$ Further note that if $\Lambda := -\lambda \log p(y)$, we obtain Proposition 3 from \cite{HASPLID}. By the exponential characterization of HASPLID and its record values, we invoke results of both (\cite{Records}, \cite{nagaraja_nevzorov_1997}) which establish the independence of the process $(Y_{R(k)} - Y_{R(k+1)})_{k \geq 0}$ from the sequence of inter-record times $(R(k+1) - R(k))_{k \geq 1}$. \\
\indent Using this fact, we introduced the new stochastic process (3.1) which we call the \textit{record improvement slope process} or \textit{slope process} for brevity. Its name comes from its expression as the \textit{ratio} of two HASPLID-specific quantities that characterize vertical, record value improvement, and horizontal, inter-record times, movement. Geometrically one can think of this random quantity as the ``derivative" of the overall minimization evaluated tangent to a particular record. 
\indent Note that the numerator $Y_{R(k)} - Y_{R(k+1)}$ denotes the amount of improvement between record values and is positive and continuous for each $k$. The denominator $ R(k+1) - R(k) $ denotes the number of inner search iterations between two consecutive records and is a discrete random variable which only takes on values in $\mathbb{N}$ for each $k$. \\
   
   \indent In order to define the new inner termination metric which is used in tandem with DMSS' original metric, $n_{\text{RECORD}}$, we first state a result that aids in its construction. Namely, we would like to know the distribution of the slope process: \\

    \noindent \textbf{Thm 1.} The common distribution of the slope process is:  \begin{multline*}
        \mathbb{P}(S_k \leq s) = -\lambda \int_{-\infty}^\infty (p(y))^{\alpha + \lambda - 1}(-\lambda \ln p(y))^{k-1}\dv{y}p(y) \\
        \cdot \displaystyle\sum_{\delta_r \in \mathbb{N}} G(k+1, -\lambda \ln p(s\delta_r - y))(1 - (p(y))^\alpha)^{\delta_r - 1} ~ dy.
    \end{multline*}
    
\textit{Proof.}  
\begin{align}
    \mathbb{P}(S_k \leq s) ={}& \int_{-\infty}^\infty \mathbb{P}(S_k \leq s | Y_{R(k)} = y) \cdot \mathbb{P}(Y_{R(k)} = y) ~ dy \nonumber \\
    ={}& \int_{-\infty}^\infty \mathbb{P}(Y_{R(k)} - Y_{R(k+1)} \leq s \cdot (R(k+1) - R(k)) | Y_{R(k)} = y) \cdot \mathbb{P}(Y_{R(k)} = y) ~ dy \nonumber \\
     \begin{split}
         ={}& \int_{-\infty}^\infty \mathbb{P}(Y_{R(k)} - Y_{R(k+1)} \leq s \cdot (R(k+1) - R(k)) | R(k+1) - R(k) = \delta_r, Y_{R(k)} = y)  \\ 
         & ~ ~ ~ ~ \cdot \mathbb{P}(R(k+1) - R(k) = \delta_r | Y_{R(k)} = y) \cdot \mathbb{P}(Y_{R(k)} = y) ~ dy \nonumber
     \end{split}\\
     \nonumber
\end{align}
We calculate each probability component-wise. The first component is:
\begin{align}
&\mathbb{P}(Y_{R(k)} - Y_{R(k+1)} \leq s \cdot (R(k+1) - R(k)) | R(k+1) - R(k) = \delta_r, Y_{R(k)} = y )  \nonumber \\
    &= \mathbb{P}(-Y_{R(k+1)} \leq s\delta_r - y) \nonumber \\
    &= \mathbb{P}(Y_{R(k+1)} > s\delta_r - y) \nonumber \\
    &= G(k+1, -\lambda \log p(s\delta_r - y)). \nonumber \\ \nonumber \medskip \\ \nonumber
\text{The second component is:} \\ \nonumber
    & \mathbb{P}(R(k+1) - R(k) = \delta_r | Y_{R(k)} = y) \\ \nonumber
    &= (p(y))^\alpha(1-p(y)^\alpha)^{\delta_r - 1}. \nonumber
    \\ \nonumber \medskip \\ \nonumber
\text{The third component is:} \\ \nonumber
    & \mathbb{P}(Y_{R(k)} = y)  \\ \nonumber
    &= -\dv{y} G(k, -\lambda \log p(y)) \nonumber \\
    &= (-\lambda \log p(y))^{k - 1} e^{\lambda \log p(y)} \cdot \dv{y} (-\lambda \log p(y)) \nonumber \\
    &= (-\lambda \log p(y))^{k-1} \cdot (p(y))^\lambda \cdot \left(-\frac{\lambda}{p(y)} \right) \cdot \dv{y}p(y). \nonumber \\
    \nonumber
\end{align}

\noindent Substituting these quantities into the above integrand and simplifying gives the desired result. 
$\square$ \\

\indent Providing a closed-form solution to the distribution above may not, and is likely not, possible. Thus, for algorithmic purposes, we turned to the expectation of the slope process. In order to levy our intuition on the metric, we condition the expected slope by the previous record value. 

\medskip

\noindent \textbf{Thm 2.} The conditional expectation of the $k$th slope value is given by:  
\begin{equation}
    \mathbb{E}[S_k | Y_{R(k)} = y] = \alpha \cdot \frac{(p(y))^\alpha}{\lambda}.
\end{equation}
\textit{Proof.}  
\begin{align}
    \mathbb{E}[S_k | Y_{R(k)} = y] &= \mathbb{E}\left[ \frac{Y_{R(k)} - Y_{R(k+1)}}{R(k+1) - R(k)} \bigg| Y_{R(k)} = y \right] \nonumber \\
    \nonumber 
\end{align}
\noindent By the memory-less property of the exponential record value differences and the continuity of their common distribution, we invoke independence of $Y_{R(k)}-Y_{R(k+1)}$ and $R(k+1) - R(k)$ which allows us to use the Law of the Unconscious Statistician:
    \begin{align*}
\mathbb{E}[S_k | Y_{R(k)} = y] =& \mathbb{E}\left[ \frac{Y_{R(k)} - Y_{R(k+1)}}{R(k+1) - R(k)} \bigg| Y_{R(k)} = y \right] \nonumber \\
    =& \mathbb{E}[Y_{R(k)} - Y_{R(k+1)} | Y_{R(k)} = y] \cdot \mathbb{E}[\left( R(k+1) - R(k) \right)^{-1} | Y_{R(k)} = y] \nonumber \\
    =& \displaystyle \frac{1}{-\lambda \log p(y)} \cdot \sum_{\delta_r \in \mathbb{N}} \frac{1}{\delta_r} \cdot \mathbb{P}(R(k+1) - R(k) = \delta_r | Y_{R(k)} = y)\nonumber \\
    =& \displaystyle \frac{1}{-\lambda \log p(y)} \cdot \sum_{\delta_r \in \mathbb{N}} \frac{1}{\delta_r} \cdot (p(y))^\alpha(1-(p(y))^\alpha)^{\delta_r - 1} \nonumber \\
    =& \displaystyle \frac{(p(y))^\alpha}{-\lambda \log p(y)} \cdot \sum_{\delta_r \in \mathbb{N}} \frac{1}{\delta_r} \cdot (1-(p(y))^\alpha)^{\delta_r - 1} \nonumber \\ 
    =& - \frac{(p(y))^\alpha}{\lambda \log p(y)}\cdot -\log\left((p(y))^\alpha\right) \nonumber \\
        =& ~ \frac{\alpha \cdot (p(y))^\alpha}{\lambda} ~ ~ \square \\
\end{align*}

\indent This is exactly the metric that we impose as an inner terminating condition on top of the already existing $n_{\text{RECORD}}$ metric. For the $k$th inner search that results in a record, we update the slope process as defined in (3.1) and, like $n_{\text{RECORD}}$, compare it to its expectation. Specifically, we bound the actual slope by its expectation. This serves two purposes, the first is that it ensures the algorithm's computational budget is not being wasted sitting in ``plateaus". Even if the $n_{\text{RECORD}}$ value has not been exceeded by the number of raw iterates, the small difference in magnitude of slopes will force $S_k \approx 0$ and thus, terminate the inner loop. The second is that it allows us to record the average ``rate" of improvement as the inner loop progresses which we can use as an analytical measure to potentially exploit in a later run.  \\

\indent Before incorporating the expected slope metric, which we denote $\mathbb{E}[S_k|Y_{R(k)}]$ for the sake of convenience, we must address several computational considerations. The first is that of $\alpha/\lambda$ that appears in (3.2). Luckily, since the outer search was unchanged, we can simply approximate this by $\zeta^{-1}$ as the value of $\zeta$ continues to update per outer search by (2.2). The second computational issue is the estimation of $p(y)$. Recall that this is defined as the CDF of the sampling measure $\rho$ which is a parameter of the conceptual algorithm HASPLID. The primary motivation of DMSS, which is \textit{modeled by} HASPLID, was to avoid computing the measure $\rho$ which is almost never available in practice, nor would its estimation be worth the computational expense. \\

\indent Unfortunately, since the expected slope includes the magnitude in the difference of record \textit{values} along with record times, as opposed to solely record times, we cannot invoke Proposition 5 to vanish $p(y)$. Instead, we make some observations to determine a suitable approximation for it. The first observation is that $p(y)$ is a CDF; by definition, it must admit a value between $0$ and $1$. The second is that $\alpha \in [0,1]$ and $\lambda \in \mathbb{R}^+$ means that the overall expectation is ``well-conditioned" in the sense that $$\frac{\alpha(p(y+\tau))^\alpha}{\lambda} - \frac{\alpha(p(y))^\alpha}{\lambda} \rightarrow 0 \quad \text{as } \tau \rightarrow 0$$and for most reasonable values of $\lambda$. Note that $\lambda$ denotes the improvement distribution parameter and thus we expect it not to be trivially close to $0$. With these observations in mind, it is reasonable to approximate $p(y)$ by any exponential, logarithmic or sufficiently smooth activation function that is frequently used in machine learning literature. A particularly tame approximation is 
\begin{equation}
    \Tilde{p}(y) = 1 - e^{-y} 
\end{equation}
which is a modified sigmoid function. Note that, because \cite{HASPLID} uses a bettering probability $(p(y))^\alpha$, we expect $p(y)$ to increase as $y$ increases. In other words, the greater $Y_{R(k)} = y$ is relative to the minima of the objective, the greater the probability of improvement should be. This approximation reflects that. Another alternative is the hyperbolic tangent function. We name the functional approximation of $p(y)$ to be $\Tilde{p}$.

\indent One final note is that $S_k$ and $E[S_k|Y_{R(k)}]$ are indexed by $k$ and not $k_R$ because both are native to a particular run (inner loop). Once a restart is performed, we assign $k_R \leftarrow k$ and refresh the slope index $k$. 

\section{Inner Search (RDMSS)}

As explained in the beginning of this chapter, we fix RDMSS' inner search to be a single iteration of a deterministic, gradient-based search method. In the experiments that follow in Chapter 4 we test the Newton-Conjugate Gradient algorithm. This is to evaluate the effectiveness of the slope metric for dealing with gradient-based local search methods. \\

\begin{figure}
     \centering
     \includegraphics[width=0.85\textwidth]{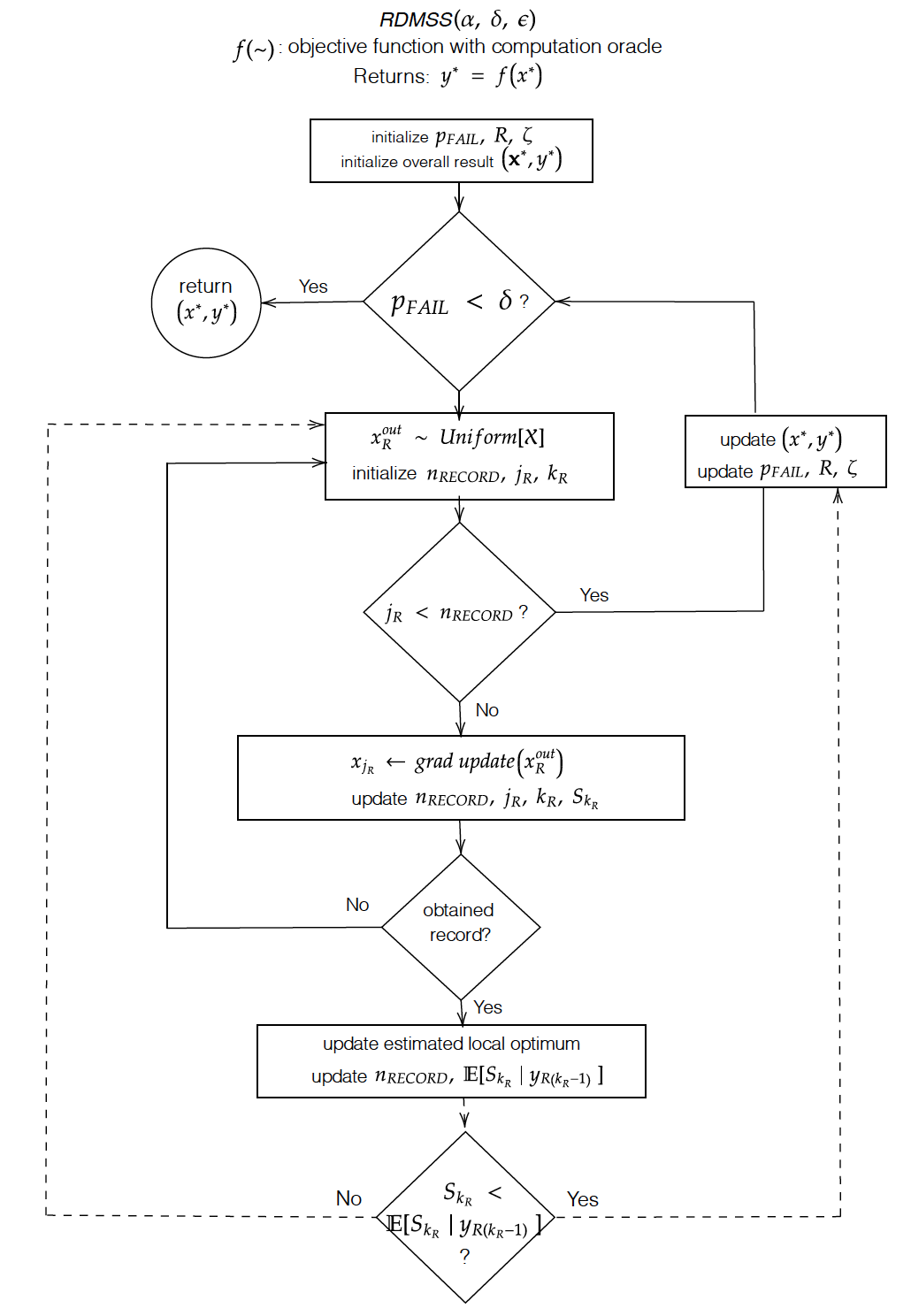}
     \caption{RDMSS($\alpha, \delta, \epsilon$). Dashed line indicates additions to DMSS.}
     \label{fig:rdmss_flow}
 \end{figure}

\indent Finally, we illustrate each major ``piece" of both DMSS and RDMSS in the following table and consolidate them into the subsequent flowchart (Figure \ref{fig:rdmss_flow}) and pseudo-code (Algorithm \ref{Algorithm 3}). \\

\begin{table}
\centering 
\begin{tabularx}{0.9\textwidth} { 
  | >{\centering\arraybackslash}X 
  | >{\centering\arraybackslash}X 
  | >{\centering\arraybackslash}X | }
  \hline
  \textbf{Multi-start ``Piece"} & \textbf{DMSS} & \textbf{RDMSS}\\
 \hline
 \hline
 Outer Termination Criteria & $p_{\text{FAIL}} < \delta$ & $p_{\text{FAIL}} < \delta$ \\
 \hline
 Outer Termination Metrics  & $p_{\text{FAIL}}, ~ R, ~ \zeta$ & $p_{\text{FAIL}}, ~ R, ~ \zeta$ \\
\hline
 \hline
 Outer Search Method & i.i.d uniform random sample over $\mathcal{X}$ & i.i.d uniform random sample over $\mathcal{X}$ \\
 \hline
  \hline
 Inner Termination Criteria & $j_R < n_{\text{RECORD}}$ & $j_R < n_{\text{RECORD}}$, and $S_k < \mathbb{E}[S_k | Y_{R(k)}]$ \\
 \hline
 Inner Termination Metrics  & $n_{\text{RECORD}}, ~ j_R, ~ k_R$ & $n_{\text{RECORD}}, ~ j_R, ~ k_R, ~ S_{k} $ and $ \mathbb{E}[S_{k} | Y_{R(k)}] $ \\ 
\hline
 \hline 
 Inner Search Method & Newton-Conjugate Gradient & Newton-Conjugate Gradient \\
\hline
\end{tabularx}\\ 
\caption{Algorithmic feature comparison table}
\label{tab:piece_wise_table}
\end{table}

 \begin{algorithm}
 \KwData{$f:$ objective function with computational oracle}
 \KwResult{$(\textbf{x}^*, \textbf{y}^*)$: global minimum $\textbf{y}^* = f(\textbf{x}^*)$ and associated location $\textbf{x}^*$}
 \noindent\rule{12cm}{0.4pt}\\ 
\textbf{Algorithm}($f$): \\
initalize $p_{\text{FAIL}}, R, \zeta$ \\
initialize $(\textbf{x}^*, \textbf{y}^*)$ \\
\While{$p_{\text{FAIL}} < \delta$}{
$y_R^{\text{out}} \leftarrow Uniform[\mathcal{X}]$ \\
initialize $n_{\text{RECORD}}, j_R, k_R, S_{k_R}$ \\
\While{$j_R < n_{\text{RECORD}}$}
{
$x_{j_R}$ $\leftarrow$ grad update$(x_R^{\text{out}})$ \\
update $j_R$ \\
\If{a record was obtained}
{update estimated local optimum \\
update $k_R, n_{\text{RECORD}}, S_{k_R} $ \\
\If{$S_{k_R} < \mathbb{E}[S_{k_R}|Y_{R(k_R-1)}]$}{break}
}
}
update $(\textbf{x}^*, \textbf{y}^*)$ \\
update $p_{\text{FAIL}}, R, \zeta$
}
\Return $(\textbf{x}^*, \textbf{y}^*)$ \\
\caption{RDMSS}
\label{Algorithm 3}
\end{algorithm}

\section{Numerical Results}

In this chapter we subject Revised DMSS to a series of tests designed to address motivating questions. Our primary task is to \textit{compare RDMSS and DMSS when equipped with a deterministic, gradient-based inner search}, specifically, the Newton-Conjugate Gradient method. This was original reason why we modified DMSS and so we prioritize the first set of experiments to answer that question. The second set of experiments is designed to compare RDMSS as a stand-alone framework against its local search method with \textit{no restarts}. For example, if we fixed the Newton-Conjugate Gradient method as the inner search for RDMSS, we would like to see how RDMSS compares to a full run of the Newton-Conjugate Gradient method with its native optimality condition. The purpose of this test is to observe a baseline comparison between RDMSS, as a representative of the multi-start framework, against a purely deterministic optimization method. Lastly, we aim to deduce the effects of dimensionality on RDMSS for the Newton Conjugate-Gradient inner search method. 

\section{Performance: Multi-start Comparison (RDMSS vs. DMSS)}

To tackle the primary objective, we first fix the inner searches of both DMSS and RDMSS to be the Newton-Conjugate Gradient (NC-G) iterative unconstrained nonlinear optimization method. RDMSS uses both $n_{\text{RECORD}}$ and $\mathbb{E}[S_{k_R} | Y_{R(k_R)}]$ as its inner loop termination metrics while DMSS only uses $n_{\text{RECORD}}$. \\
\indent We perform our experiments on several objective test functions where $f: \mathcal{X} \subset \mathbb{R}^5 \rightarrow \mathbb{R}$. We fix $\mathcal{X}$ to be the subset of 5-dimensional Euclidean space upon which our oracle is defined. We also fix both algorithm's parameters as $\epsilon = (0.01)^5, \delta = 0.001, $ and $\alpha = 0.5$. For each objective test function, we plot 50 individual \textit{function histories} of both DMSS and RDMSS, distinguished by color. The function history is simply a record of every time the algorithm used the computational oracle $f(.)$ and what its output was. In other words, we are plotting every raw search iterate across every inner loop until the algorithm is globally terminated. The objective test functions we experimented were all standard test functions with known global optima. They are the: Zakharov, Rosenbrock, Rotated Hyper-Ellipsoid, Styblinski-Tang, Shifted-Sinusoidal and Centered Sinusoidal test functions, as specified in \textit{Appendix A}. \\
\indent For each test function, we observe the overall error between the returned solutions of DMSS/RDMSS and the known global optimum. We also include several metrics which are typically used to characterize the performance of multi-start algorithms on any given \textit{global} run. These are: 
\begin{enumerate}
    \item Total number of restarts 
    \item Number of function evaluations required to reach the global optimum 
    \item Average number of raw search iterations per inner loop 
    \item Total number of function evaluations 
\end{enumerate}
For the set of all fifty runs per algorithm, we aggregate the above metrics as averages and we also calculate 
\begin{enumerate}
    \item Ratio of runs that entered the $\epsilon-$target region to those that did not (success \%)
\end{enumerate} 

\noindent These metrics are consolidated in Table 4.1.

\subsection{Zakharov}
\begin{figure}[H]
    \centering
    \includegraphics[width=\textwidth]{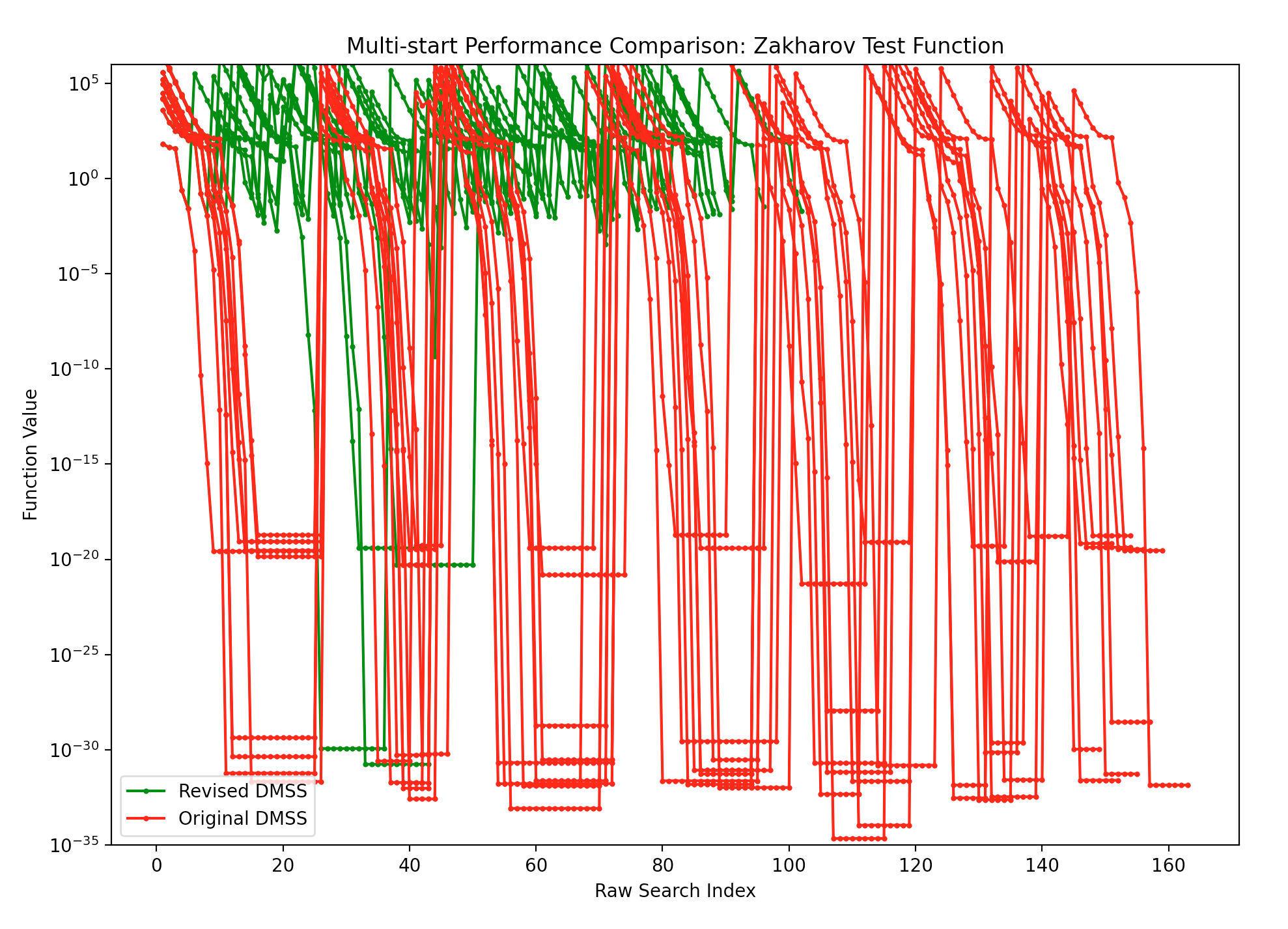}
    \caption{Raw function histories for RDMSS (green) and DMSS (red) applied to the Zakharov function}
    \label{fig:zakharov_performance}
\end{figure}

In Figure \ref{fig:zakharov_performance} we notice that RDMSS is successful in reducing the number of overall function evaluations and increasing the frequency of restarts. Additionally, it does not lose accuracy. All of the RDMSS runs reach the global optimum as supported by the success rates given in Table 4.1. The plots in Figure 4.1 clearly show that RDMSS restarts much earlier than DMSS when performing inner searches. In other words, the slope criterion makes the algorithm perform a restart much sooner, only giving each inner loop just enough time to sufficiently explore its particular local search neighborhood and find an optimum. We hypothesize that for other test functions, we will continue seeing a sharp reduction in the average number of inner loop iterations. 

\subsection{Rosenbrock}
\begin{figure}[H]
    \centering
    \includegraphics[width=\textwidth]{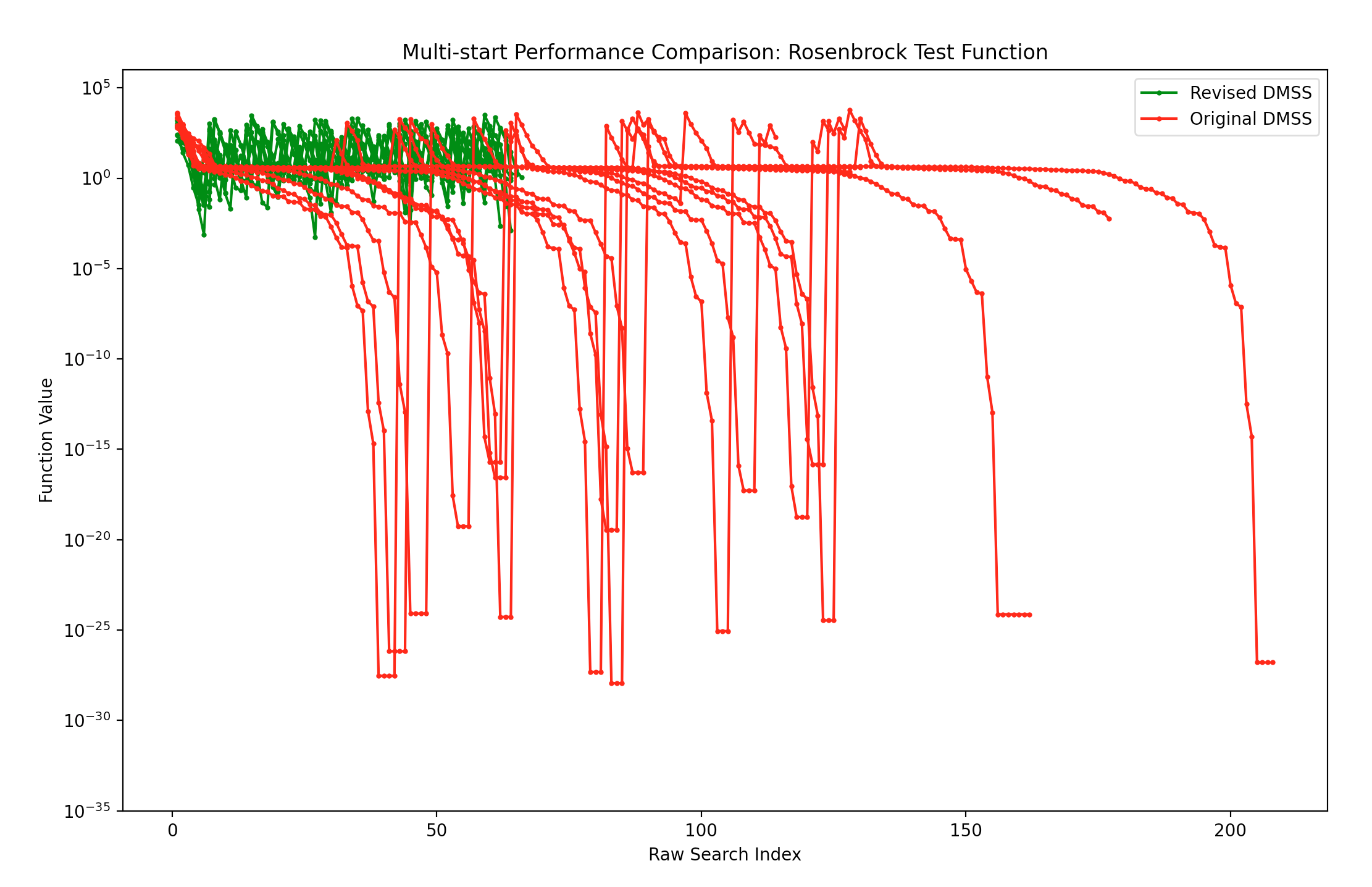}
    \caption{Raw function histories for RDMSS (green) and DMSS (red) applied to the Rosenbrock function}
    \label{fig:rosen_performance}
\end{figure}

\begin{figure}[H]
    \centering
    \includegraphics[width=\textwidth]{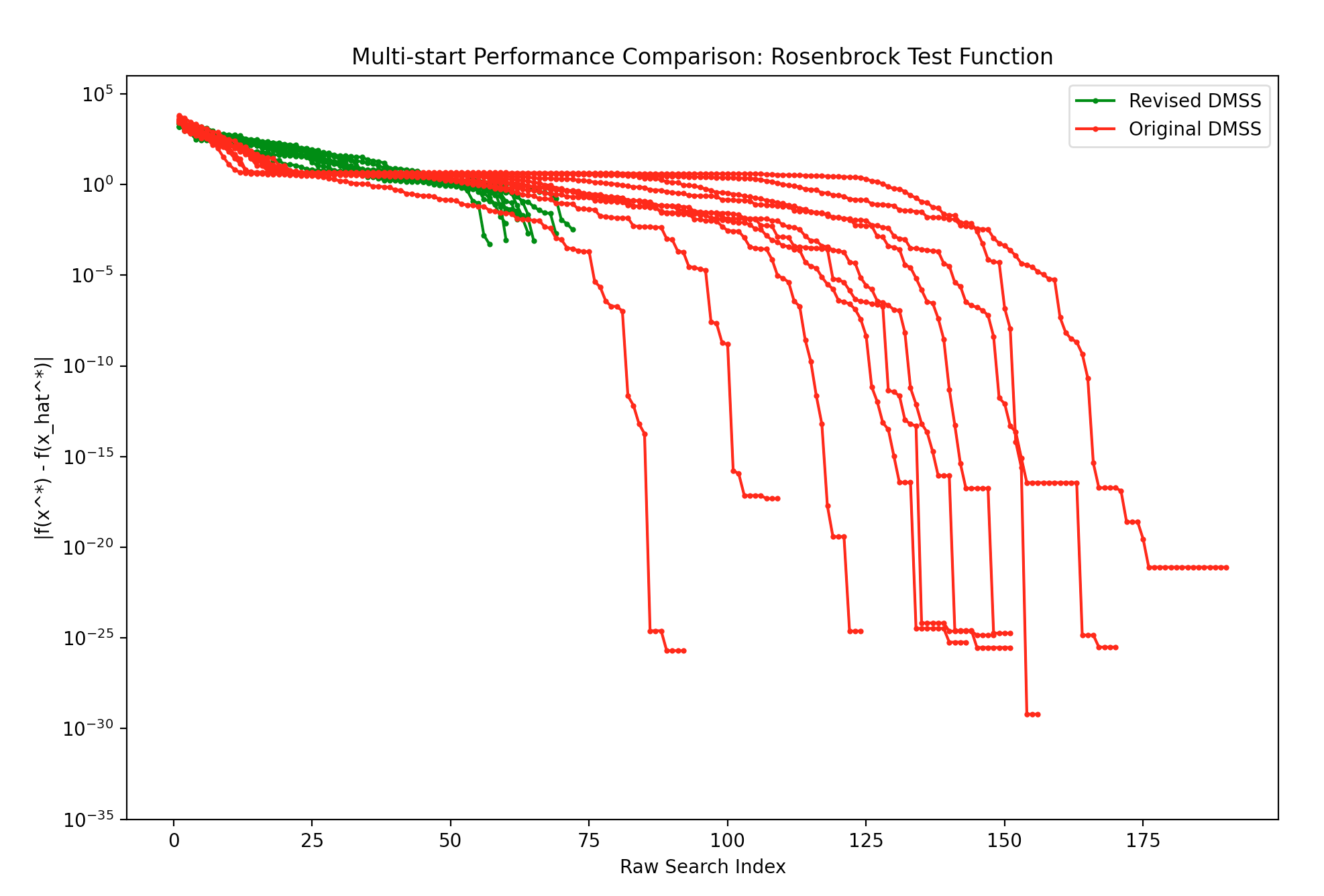}
    \caption{Sorted function errors for RDMSS (green) and DMSS (red) applied to the Rosenbrock function}
    \label{fig:rosenbrock_sorted}
\end{figure}

The Rosenbrock test function, whose experiment is given in Figure \ref{fig:rosen_performance}, does feature a reduction in the total number of function evaluations, but its average number of inner searches per inner loop reduces too dramatically, at the cost of success rate. However, this is expected behavior. The Rosenbrock function is ``valley" shaped while the Zakharov function is ``plate-shaped". The initially ``steep" descent to the global minimum at the origin followed by an increasingly ``flat" hyper-dimensional surface means that initially biased expected slope will admit a soft inflection point. It is likely that the slope measurement is prematurely considered ``bad" by the inner search criterion and then activated as a result. To verify this, we take the same runs given in Figure \ref{fig:rosen_performance} and plot their sorted errors instead of their raw function histories, see Figure \ref{fig:rosenbrock_sorted}. Even though the plot depicts a case of RDMSS being ``too impatient", we note that the target region for five dimensions is given by $[-(0.01)^5, (0.01)^5] = [- 10^{-10}, 10^{10}]$. The RDMSS runs reach an estimated optimum between $10^{-3}$ and $10^{-5}$ (estimated visually), the difference in true optimum $f(x^*) = 0$ and the estimated optimum $f(\hat{x}^*) \approx 10^{-5}$ against the number of evaluations required to reach an improvement is negligible i.e., corresponds to a near-zero slope and thus, decides to terminate. It is clearer that objectives with ``flattening" around a global optimum will exhibit similar behavior. 

\subsection{Rotated Hyper-Ellipsoid}
\begin{figure}[H]
    \centering
    \includegraphics[width=\textwidth]{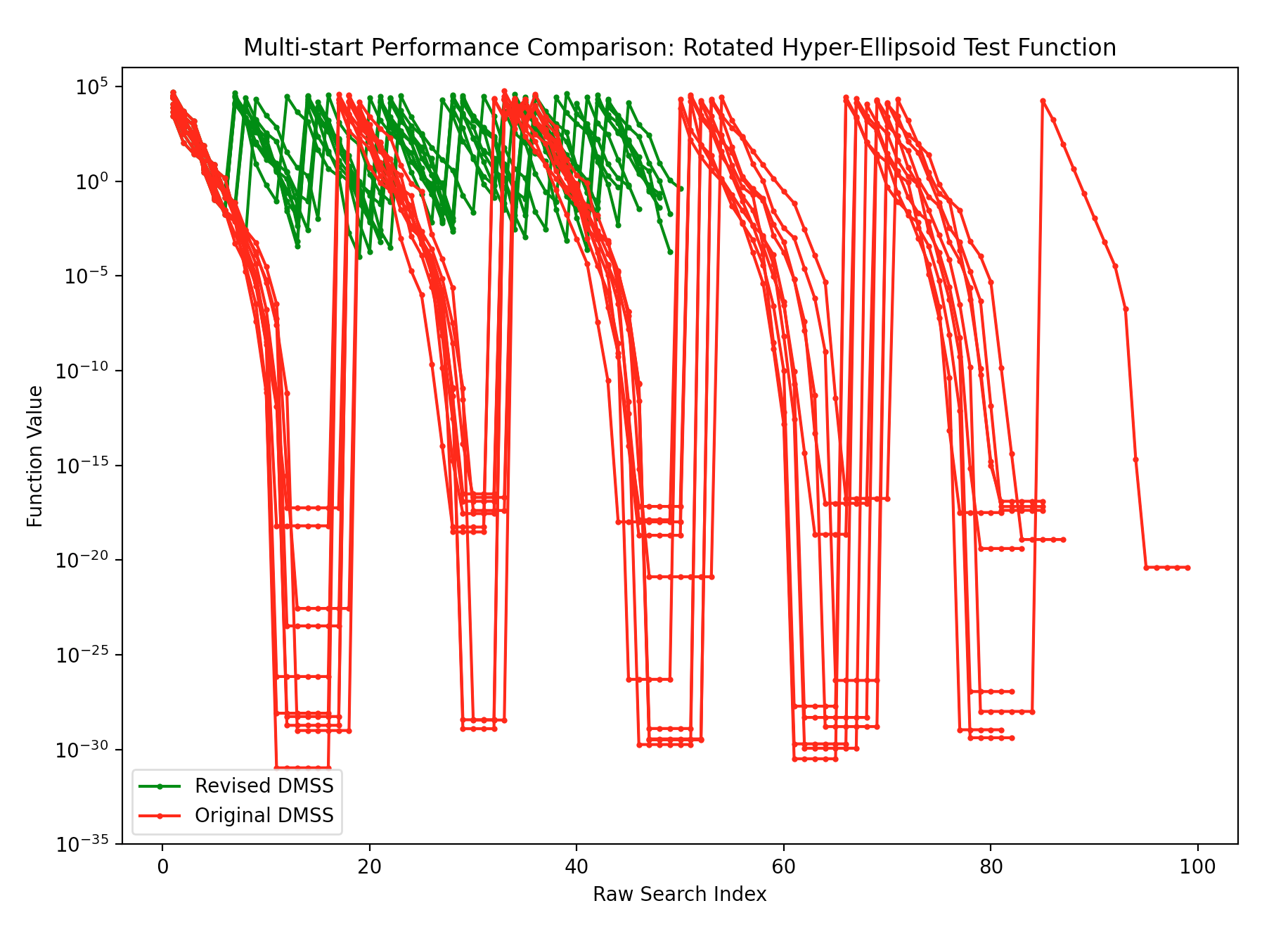}
    \caption{Raw function histories for RDMSS (green) and DMSS (red) applied to the Rotated Hyper-Ellipsoid function}
    \label{fig:RHE_performance}
\end{figure}

The Rotated Hyper-Ellipsoid function in Figure \ref{fig:RHE_performance} is another example of RDMSS being too aggressive in its inner search termination criteria. Like the Rosenbrock function, the Rotated Hyper-Ellipsoid surface begins with a sharp drop-off and levels out as it reaches its global minimum at the origin. In fact the Rotated Hyper-Ellipsoid function, also known as the sum of squares function, is convex and defined on a relatively large domain $[-65.536, 65.536]$, which once again supports the geometric argument presented in the previous analysis of the Rosenbrock experiments. As in the previous experiments, we see a too sharp of a decrease in average number of inner search iterations, leading to a 0 success rate, given in Table 4.1.

\subsection{Styblinski-Tang}
\begin{figure}[H]
    \centering
    \includegraphics[width=\textwidth]{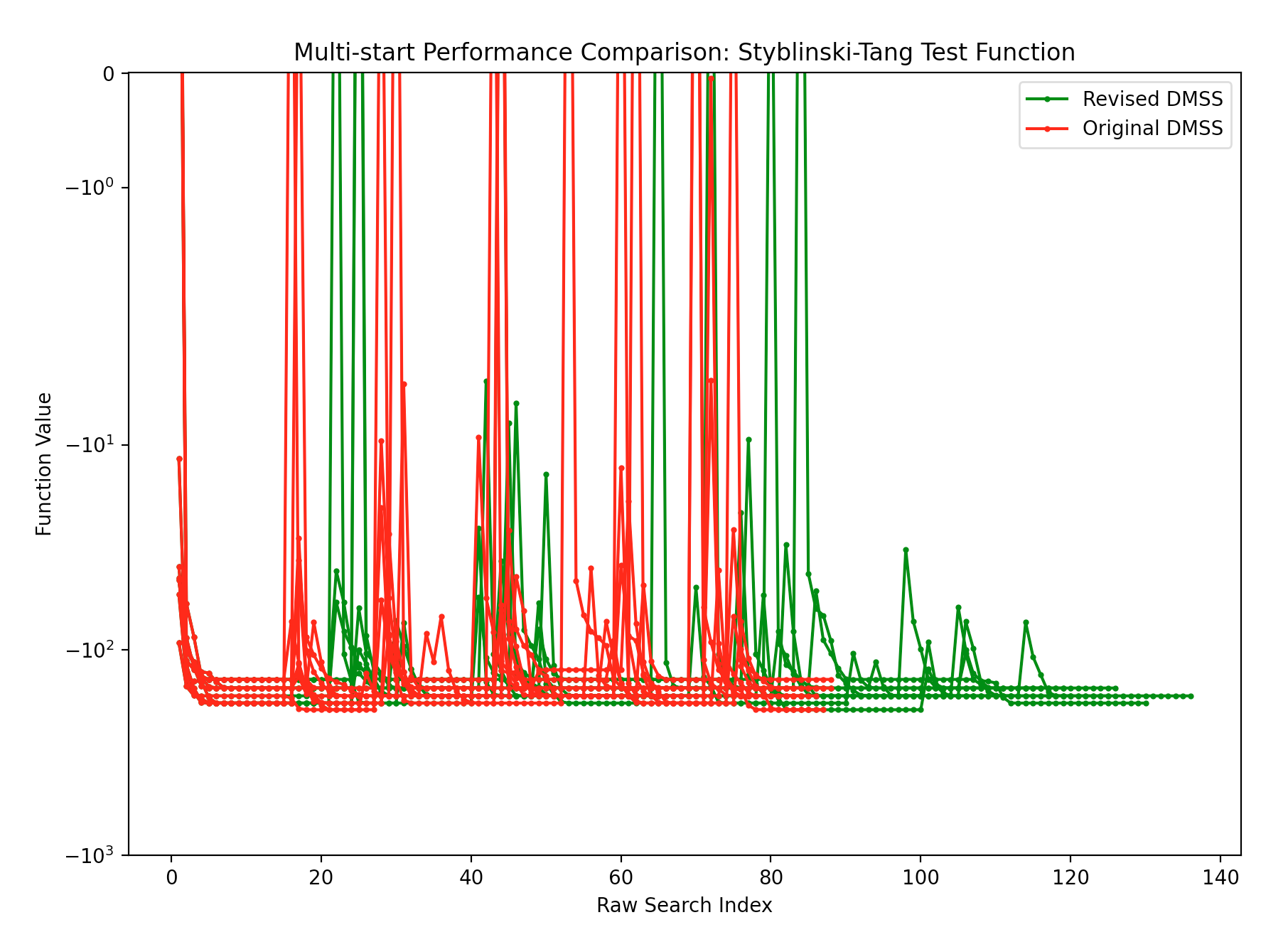}
    \caption{Raw function histories for RDMSS (green) and DMSS (red) applied to the Styblinski-Tang function}
    \label{fig:styblinski_tang_performance}
\end{figure}

The Styblinski-Tang experiment shown in Figure \ref{fig:styblinski_tang_performance} deviates from the previous three test functions. We see that RDMSS actually takes ``longer" to terminate globally and Table 4.1 shows that RDMSS actually reduced the total number of restarts by nearly 1. Again, shape likely plays a role here but herein lies the power of multi-start frameworks in general. Unlike the previous three test functions, Styblinski-Tang is multimodal. Its overall shape is indeed roughly bowl-shaped however the presence of local minima means that there are several inflection points which were circumvented by restarts. This objective highlights the advantage of the slope metric as opposed to the $n_{\text{RECORD}}$ metric as record-difference or ``jump" (\textit{see Appendix B}) magnitude becomes an important variable by which the algorithm decides whether or not to traverse the domain, i.e., perform a restart. 

\subsection{Shifted Sinusoidal}
\begin{figure}[H]
    \centering
    \includegraphics[width=\textwidth]{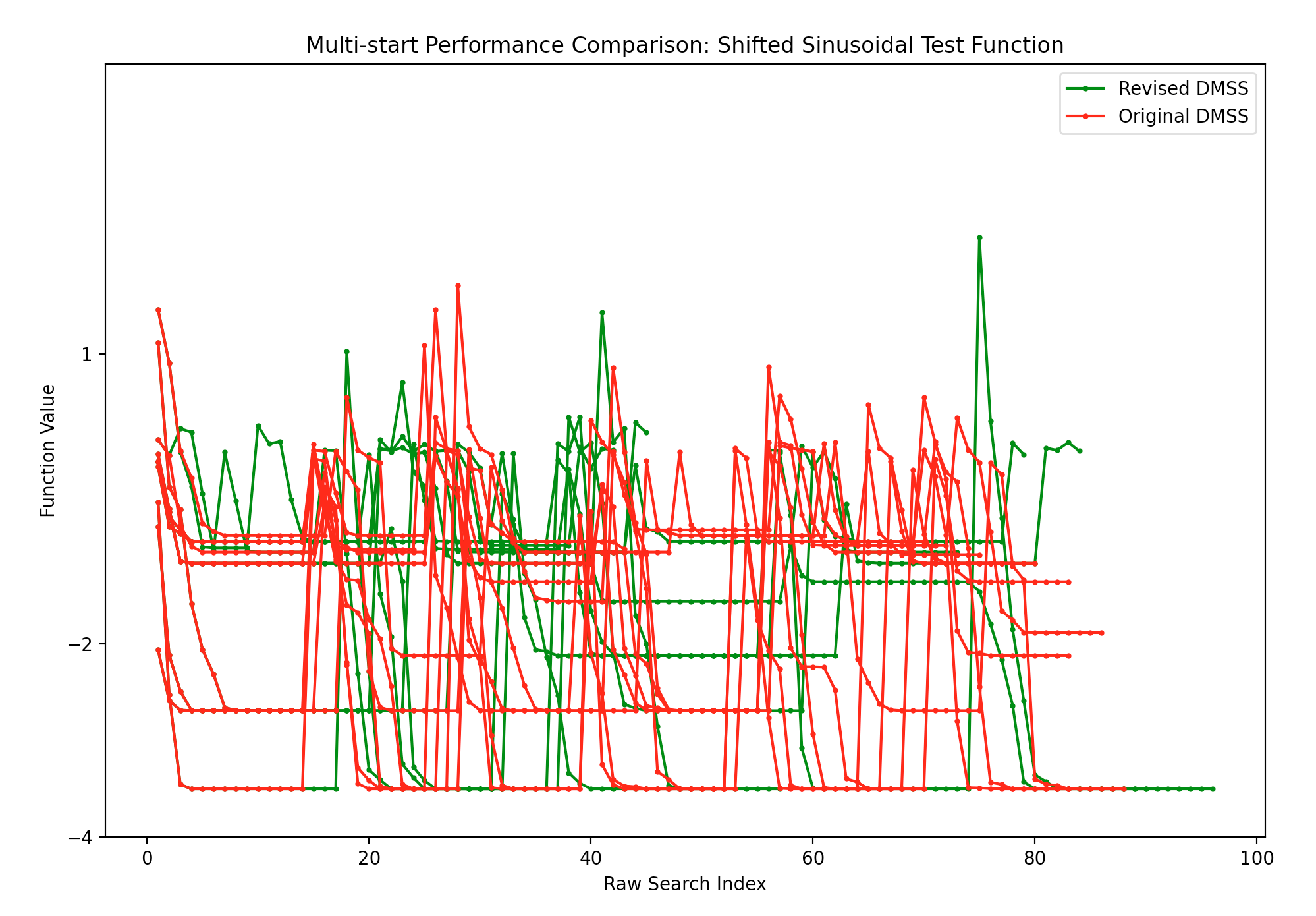}
    \caption{Raw function histories for RDMSS (green) and DMSS (red) applied to the Shifted Sinusoidal function}
    \label{fig:shifted_sine_performance}
\end{figure}

The Shifted Sinusoidal experiments also yield promising results. It is not evident from Figure \ref{fig:shifted_sine_performance}, but Table 4.1 shows that with a maintained average of 6 restarts, RDMSS was able to arrive at the $\epsilon-$target region in fewer iterations \textit{and} the average number of inner search iterations was reduced along with the total number of function evaluations required for termination. The success rates were identical, however. Once again, the shifted sinusoidal function is very multimodal. It is at this time that we can start to formulate a conclusive hypothesis that RDMSS performs better in multimodal settings since the slope criterion, which is directly correlated to the gradient step-size, takes algorithmic precedence over the number of iterations required to terminate (i.e., $n_{\text{RECORD}}$) allowing the algorithm to spend more time exploring the search domain and less time trapped in local search spaces. 

\subsection{Centered Sinusoidal}
\begin{figure}[H]
    \centering
    \includegraphics[width=\textwidth]{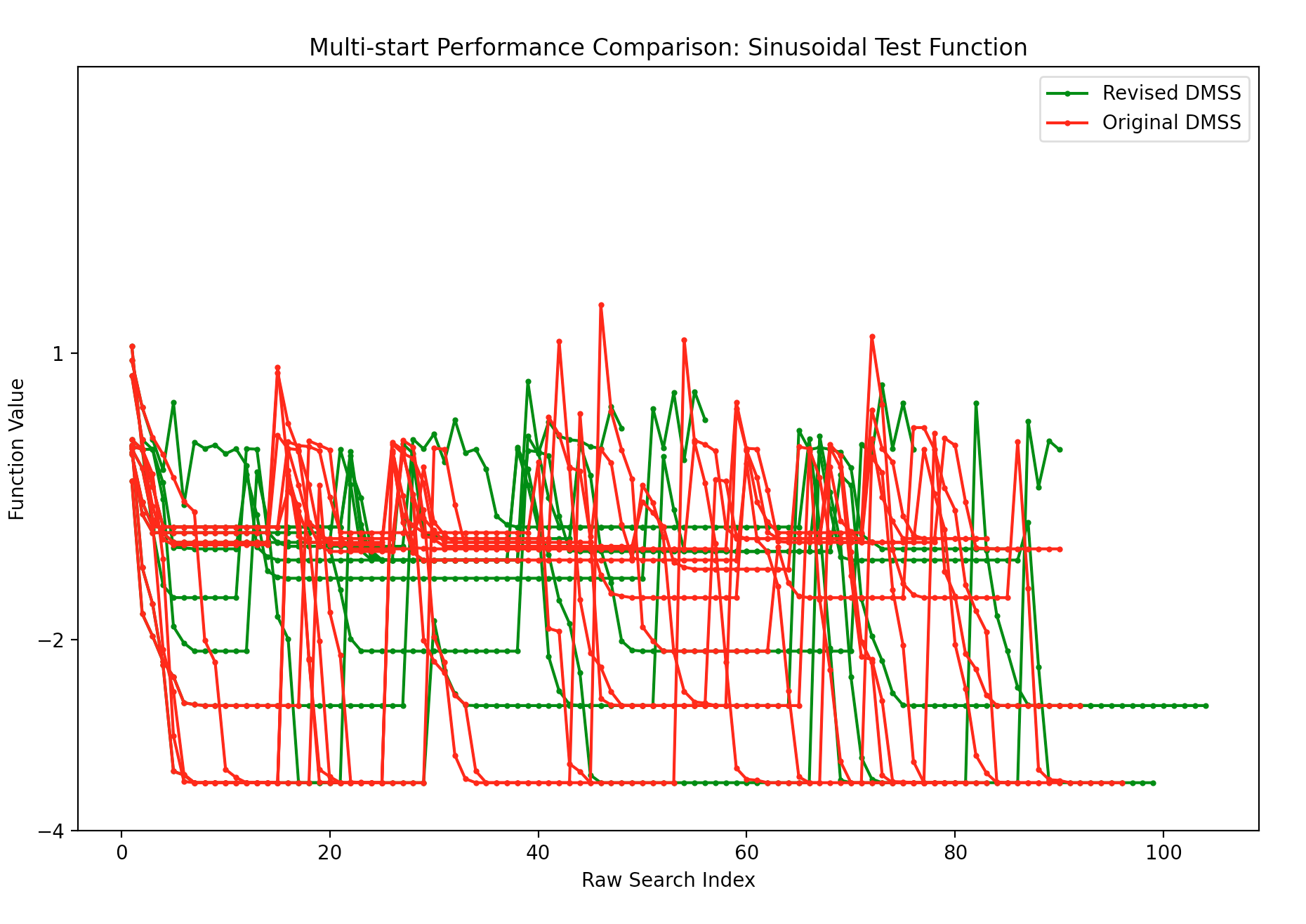}
    \caption{Raw function histories for RDMSS (green) and DMSS (red) applied to the Sinusoidal function}
    \label{fig:sine_performance}
\end{figure}

To reinforce the idea that the sharper inner termination criteria of RDMSS is advantageous for multimodal functions, we also test the centered sinusoidal objective in Figure \ref{fig:sine_performance}. Table 4.1 shows a similar pattern for both the shifted and centered sinusoidal functions, the average number of inner loop iterations decreased and the total number of restarts remained constant, resulting in less overall function evaluations. In this set of experiments, we experienced a slight reduction in the overall success rate but, as with the Rosenbrock function, this is expected since the $\epsilon$-target region is so small and our slope metric tends to zero. \\



\begin{table}
\centering
 \begin{tabular}{||c c c c c c c||}
 \hline 
\thead{Objective}  & \thead{Method} & \thead{avg. \# \\ of restarts} & \thead{avg. \# of  \\function \\ evaluations \\ to reach \\ $\epsilon-$target \\ region} & \thead{avg. \# of \\ searches \\ per inner \\ loop} & \thead{total \# of \\ function \\ evaluations} & \thead{\# of \\ successes \\ out of \\ 50} \\%
 \hline
Zakharov & DMSS & 8.08 & 69.6 & 19.78 & 159.74 & 50 \\ 
 & RDMSS & 9.56 & 19.92 & 9.40 & 108.1 & 50 \\    [1ex] 
 \hline
Rosenbrock & DMSS & 4.10 & 91.36 & 36.91 & 150.18 & 50 \\ 
 & RDMSS & 14.90 & N/A & 4.10 & 61.08 & 0 \\    [1ex] 
 \hline
\makecell{Rotated \\ Hyper-Ellipsoid} & DMSS & 5.88 & 44.12 & 15.17 & 88.46 & 50 \\ 
 & RDMSS & 7.00 & N/A & 6.60 & 46.20 & 0 \\    [1ex] 
 \hline
Styblinski-Tang & DMSS & 5.90 & 24.43 & 15.74 & 92.44 & 0 \\ 
 & RDMSS & 5.00 & 49.00 & 24.56 & 122.78 & 13 \\   [1ex] 
 \hline
\makecell{Shifted \\ Sinusoidal} & DMSS & 5.94 & 35.69 & 14.36 & 85.02 & 39 \\ 
 & RDMSS & 6.16 & 29.28 & 12.58 & 74.14 & 39 \\   [1ex] 
 \hline
 \makecell{Centered \\ Sinusoidal} & DMSS & 5.98 & 37.28 & 14.36 & 85.8 & 43 \\ 
 & RDMSS & 6.08 & 27.91 & 12.85 & 75.66 & 34 \\  [1ex]
 \hline 
\end{tabular}
 \label{multi_start_performance}
 \caption{Figures generated for $\alpha = 0.5, \delta = 0.001$ and $\epsilon = (0.01)^5$ for compact domain $\mathcal{X} \subset \mathbb{R}^5$ aggregated over 50 individual runs per method.}
\end{table} 

\bigskip 

\begin{figure}
    \centering
    \includegraphics[width=\textwidth]{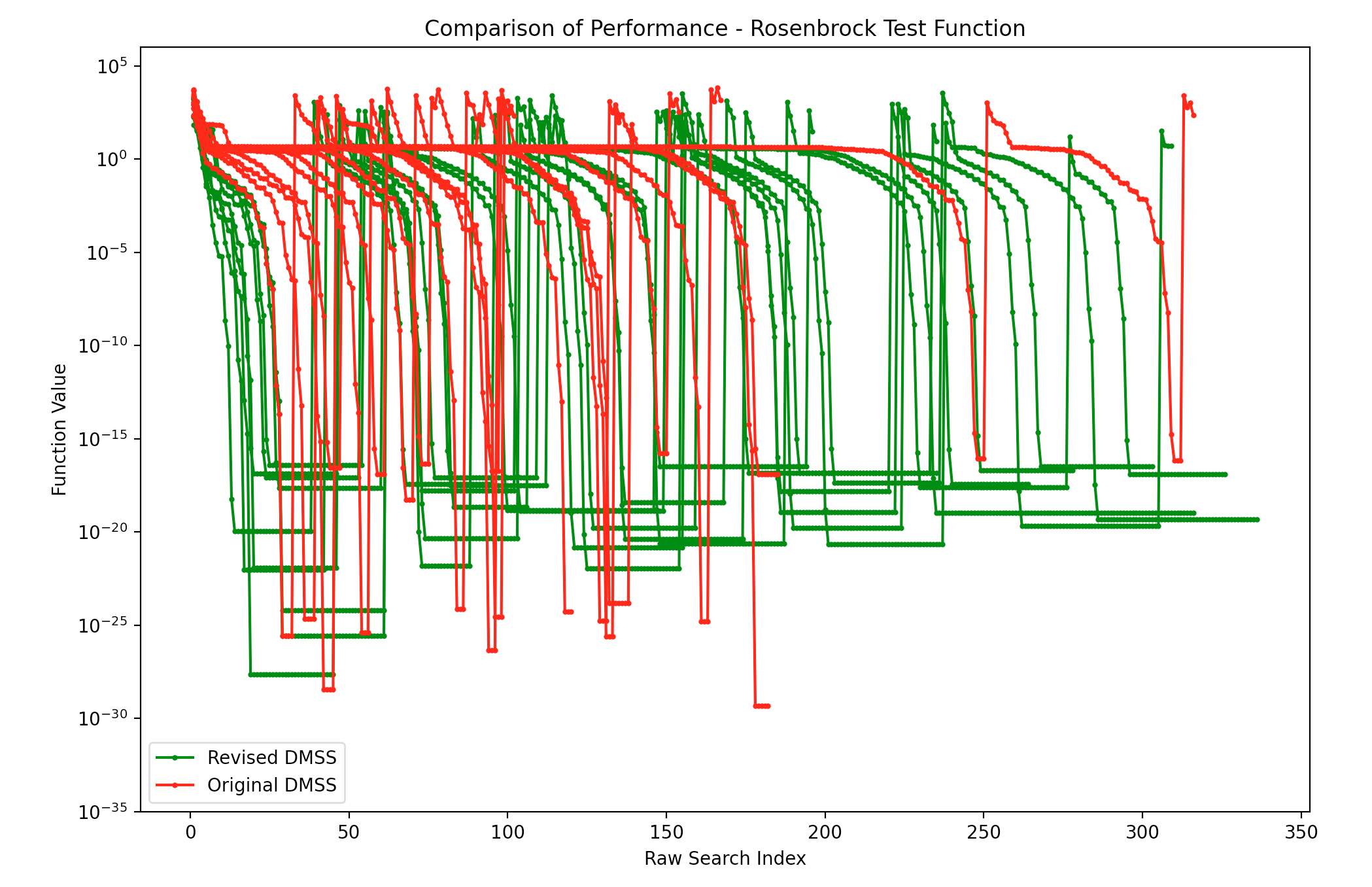}
    \caption{Raw function histories for RDMSS (green) and DMSS (red) applied to the Rosenbrock function in $\mathcal{X} \subset \mathbb{R}^5$ with $\Tilde{p}(y) = 1-e^{-y/(2^5)}$}
    \label{fig:conclusion_rosen}
\end{figure}

Before moving on to the next experiment, we would like to present a consideration on the Rosenbrock test function's numerical results which was found upon doing a deeper dive into why its results were so starkly different from the rest. Considering the form of the $\Tilde{p}$ in (3.3), the assumed functional form of the range CDF of HASPLID that models DMSS, we find that 
\begin{equation*}
    \Tilde{p}(y) = 1 - e^{-y} \approx 1
\end{equation*}
for almost any positive $y \geq 3$. When evaluating across the domain of the Rosenbrock (recall it's ``valley-shape") in high dimensions, $\Tilde{p}(y)$ is far more likely to admit a value trivially close to one, \textit{biasing} the early iterates of the slope process and subsequently cause the expected slope to be greater than it should be i.e., expect more aggressive improvements. In order to effectively ``scale" this function out in accordance to the high dimensionality of $\mathcal{X}$, we attempted the same experiment with the approximation 
\begin{equation*}
    \Tilde{p}(y) = 1 - e^{-y/(2^d)}
\end{equation*}
for $d=5$ and obtained the results in Figure \ref{fig:conclusion_rosen}. Immediately we see that RDMSS behaves almost exactly like DMSS except, strangely, it looks as if RDMSS is taking \textit{longer} than DMSS and which is exactly the opposite of what it was built for. \\
\indent However, upon further investigation, it is not strange at all. DMSS places an upper bound on the \textit{total} number of records per inner search while RDMSS places an upper bound on the number of records \textit{since the previously obtained record}. This is a subtle but important distinction because it means that the $n_{\text{RECORD}}$ metric is only applied after the ``rate" of the local search, modeled by the slope, slows down. \\
\indent In other words, by adjusting the scaling of the approximation function, we can effectively turn RDMSS into DMSS. This opens up the possibility of dynamically adjusting how ``economical" we want our algorithm to be with its computational expense when searching for local optima. With no scaling, the most conservative case is RDMSS while, with a high scaling  factor, the least conservative case is DMSS. \\

\section{Performance: Deterministic Comparison}
To address our secondary objective, we test RDMSS equipped with NC-G as its inner search against fifty unadulterated runs of NC-G with no stochastic restarts. The purpose of these tests is to highlight the importance of multi-start and perform an analysis of the trade-offs that occur in its use. Some of the factors we consider in our analysis include: computational expense, accuracy of solution and ``speed". Here we use ``speed" to mean the average number of function evaluations required to enter the $\epsilon-$target region. 

\indent The initial points for both algorithms are identical for each individual global run, this is done by means of a fixed random seed. We use the same RDMSS parameters as stated in Section 4.1 and assign our test-function domains a dimension of $5$. We run the experiment for the same set of test functions as stated in Section 4.1.

\subsection{Zakharov} 
\begin{figure}[H]
    \centering
    \includegraphics[width=\textwidth]{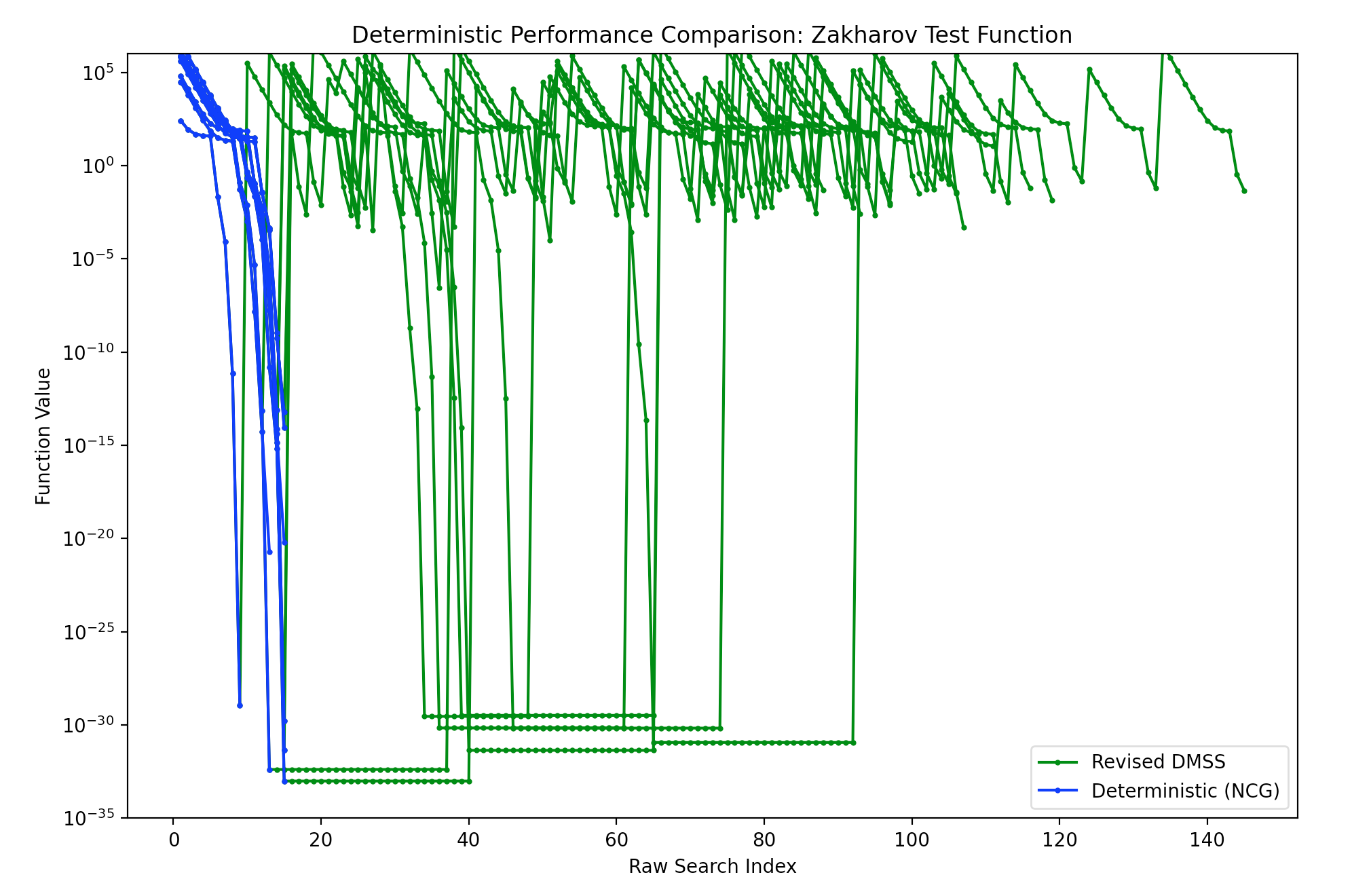}
    \caption{Raw function histories for RDMSS (green) and NCG with no restarts (blue) applied to the Zakharov function}
    \label{fig:zakharov_deterministic}
\end{figure}

The Zakharov function experiments depicted in Figure \ref{fig:zakharov_deterministic} show that the NC-G method is able to approach the global optimum very quickly and accurately. We see that on several runs, RDMSS is keeping up with NC-G and then suddenly restarts. It is likely that at this point, RDMSS decides that it is no longer worth the negligible increase in accuracy to continue searching in the current neighborhood and breaks from the inner loop to search a different local search domain. Despite the increase in number of evaluations required to reach the global optimum and the overall number of function evaluations, the success rates of RDMSS and NC-G are both 100 \%, see Table 4.2. 

\subsection{Rosenbrock}
The deterministic Rosenbrock experiments sing a very similar tune to the multi-start performance experiments. RDMSS approaches the $\epsilon-$target region but determines that the relative improvement is no longer worth the computational expense. However, as mentioned in the previous experiment, there is a way to ``scale" the termination condition for the slope via the $\Tilde{p}$. What is particularly notable about this experiment is that some of the randomly initialized starting points of the N-CG runs were poor enough that the NC-G algorithm did not converge to the $\epsilon-$target region 100 \% of the time. 

\subsection{Rotated Hyper-Ellipsoid}
The Rotated Hyper-Ellipsoid function, like the Rosenbrock function, has a shape that is not conducive to the slope criteria that RDMSS uses as its primary inner optimization mechanism. Unsurprisingly, its convexity is why NC-G has no trouble at all reaching the global optimum in very few iterations. 

\subsection{Styblinski-Tang}
The Styblinski-Tang function responded positively to this test as it did in the Multi-start performance comparison. The NC-G method, although fast, was not accurate enough and seemed to get trapped in local optimum over 90\% of the time as per Table 4.2.

\subsection{Shifted-Sinusoidal}
The Shifted Sinuosidal function, like the multimodal Styblinski-Tang function, shows in Figure \ref{fig:shifted_sine_deterministic} why multi-start is necessary. We can clearly see that the NC-G method gets trapped in local optima the majority of the time while RDMSS' restarts allow it to escape and eventually converge to the correct target region. RDMSS presents a substantial increase in accuracy, obviously at the expense of a greater number of function evaluations. 

\subsection{Centered Sinusoidal}
The Centered Sinusoidal test function exhibits behavior similar to its shifted version. Once again, we can clearly see that the NCG runs get trapped in local optima while the RDMSS runs have a much higher probability of escaping. This is supported by the metrics in Table 4.2.  

\begin{figure}[H]
    \centering
    \includegraphics[width=\textwidth]{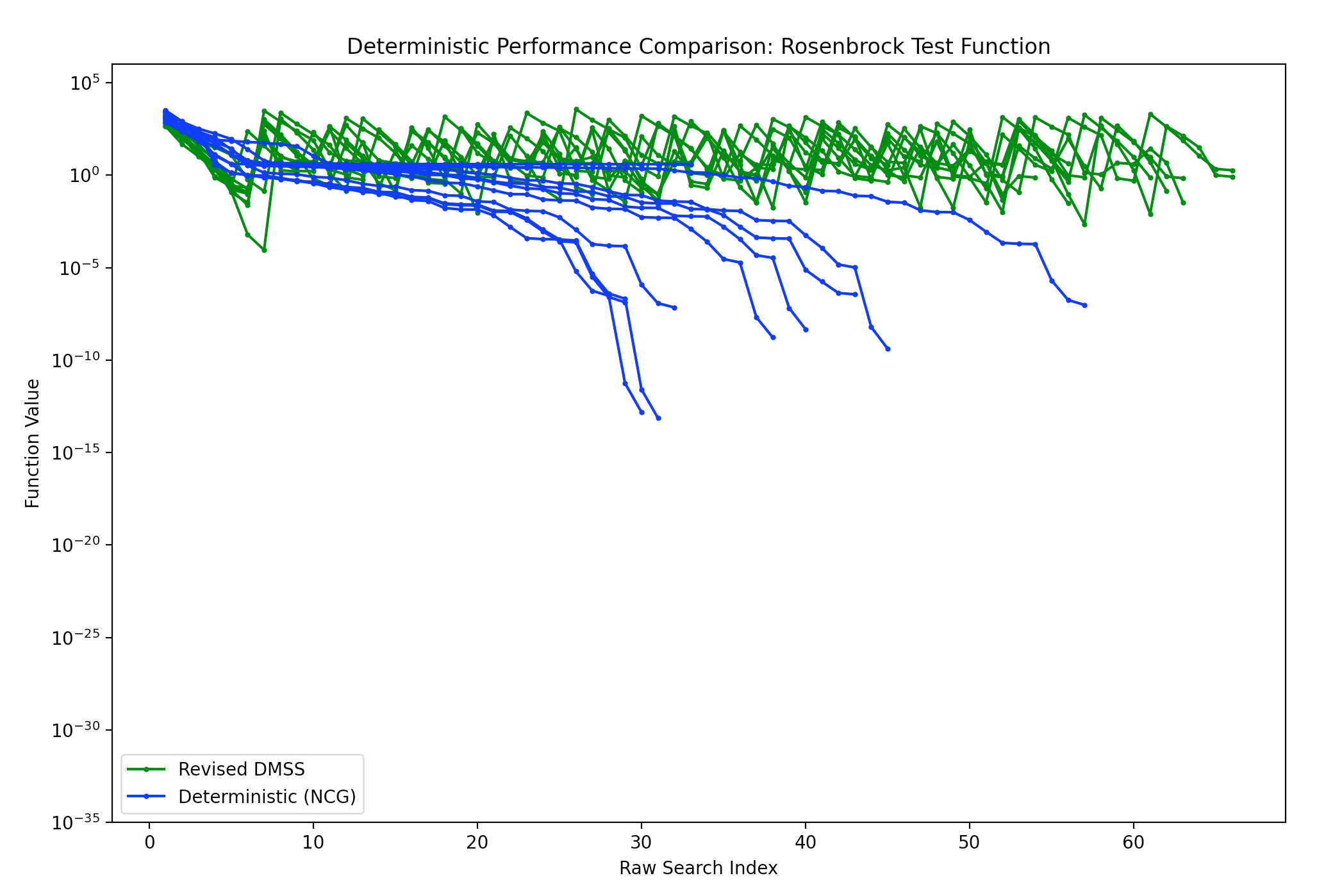}
    \caption{Raw function histories for RDMSS (green) and NCG with no restarts (blue) applied to the Rosenbrock function}
    \label{fig:Rosenbrock_deterministic}
\end{figure}

\begin{figure}[H]
    \centering
    \includegraphics[width=\textwidth]{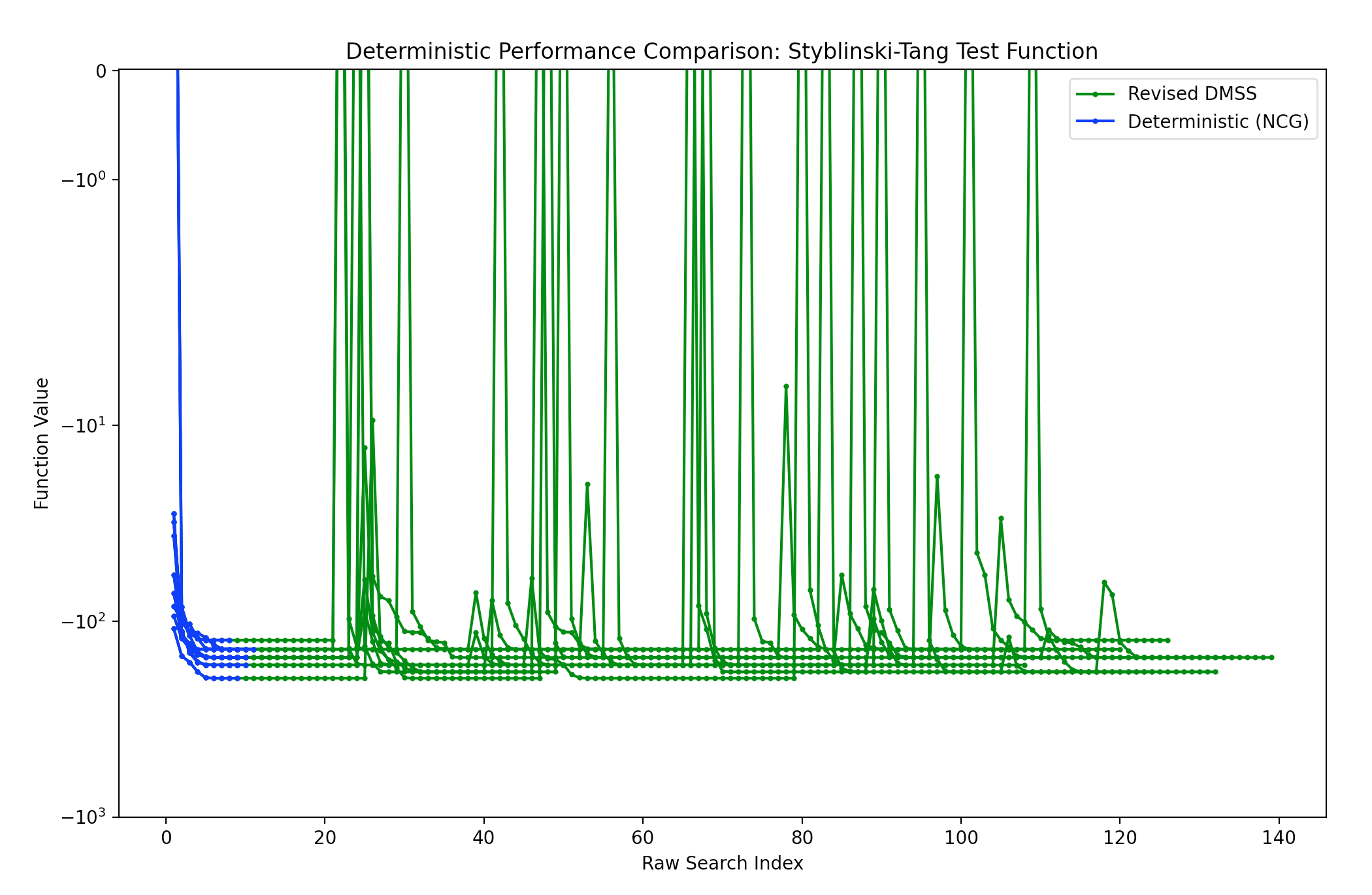}
    \caption{Raw function histories for RDMSS (green) and NCG with no restarts (blue) applied to the Styblinki-Tang function}
    \label{fig:styblinski_tang_deterministic}
\end{figure}

\begin{figure}[H]
    \centering
    \includegraphics[width=\textwidth]{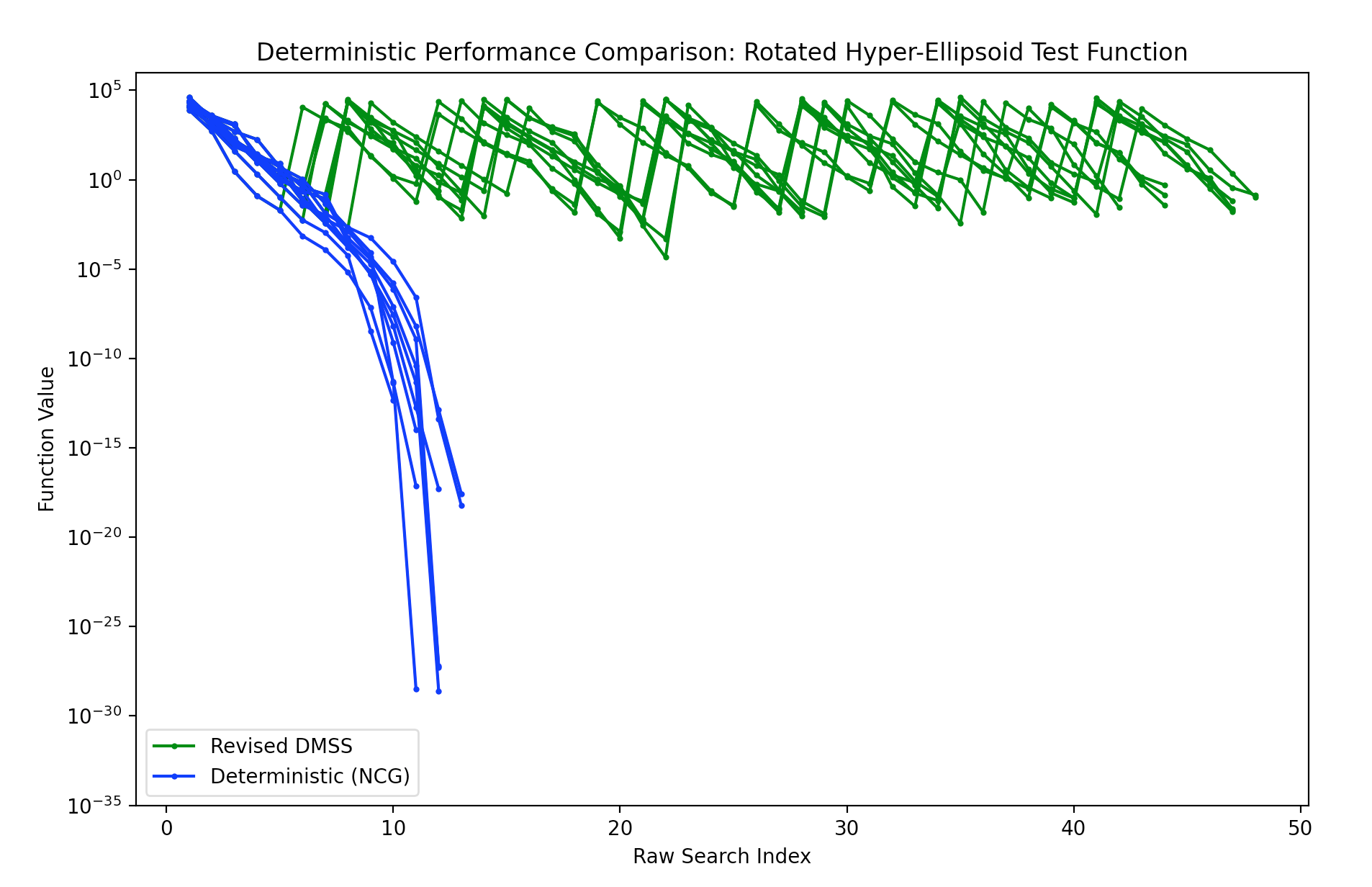}
    \caption{Raw function histories for RDMSS (green) and NCG with no restarts (blue) applied to the Rotated Hyper-Ellipsoid function}
    \label{fig:RHE_deterministic}
\end{figure}

\begin{figure}[H]
    \centering
    \includegraphics[width=\textwidth]{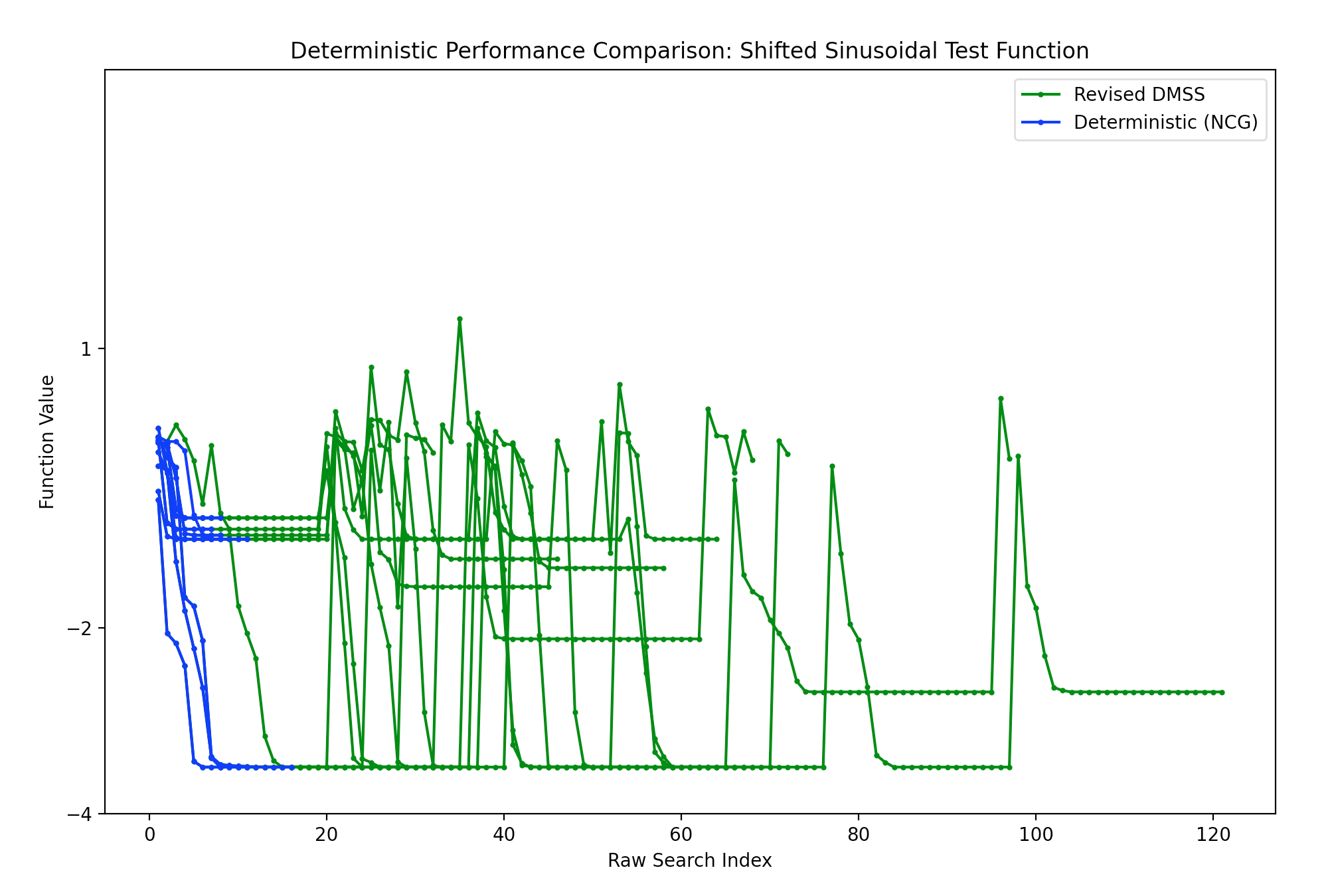}
    \caption{Raw function histories for RDMSS (green) and NCG with no restarts (blue) applied to the Shifted-Sinusoidal function}
    \label{fig:shifted_sine_deterministic}
\end{figure}

\begin{figure}[H]
    \centering
    \includegraphics[width=\textwidth]{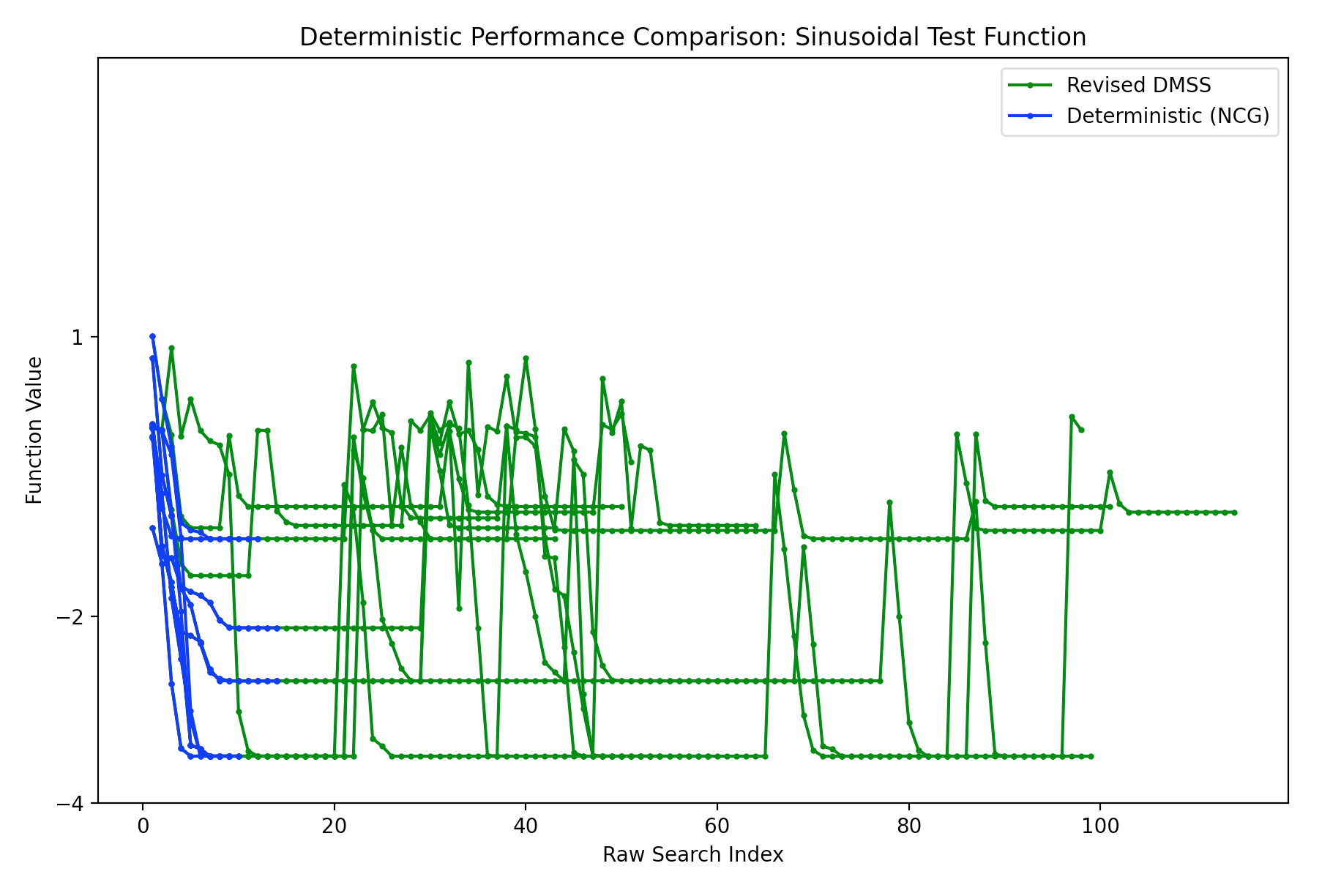}
    \caption{Raw function histories for RDMSS (green) and NCG with no restarts (blue) applied to the Sinusoidal function}
    \label{fig:shifted_sine_deterministic}
\end{figure}

\begin{table}
\centering
 \begin{tabular}{||c c c c c||} 
 \hline
 \thead{Objective} & \thead{Method} & \thead{avg. \# of  \\function \\ evaluations \\ to reach \\ $\epsilon-$target \\ region} & \thead{total \# \\ of \\ function \\ evaluations}  & \thead{\# of \\ successes \\ out of \\ 50}  \\ [0.5ex]
 \hline\hline
Zakharov & NC-G & 11.64 & 11.66 & 50 \\ 
 & RDMSS & 19.92 & 108.10 & 50 \\   [1ex] 
 \hline
Rosenbrock & NC-G & 37.33 & 47.36 & 43 \\ 
 & RDMSS & N/A & 61.08 & 0 \\   [1ex] 
 \hline
Rotated Hyper-Ellipsoid & NC-G & 10.88 & 10.90 & 50 \\ 
 & RDMSS & N/A & 46.20 & 0 \\   [1ex] 
 \hline
Styblinski-Tang & NC-G & 6.67 & 8.50 & 3 \\ 
 & RDMSS & 49.00 & 122.78 & 13 \\  [1ex] 
 \hline
Shifted Sinusoidal & NC-G & 7.33 & 8.24 & 12 \\ 
 & RDMSS & 29.28 & 74.14 & 39 \\    [1ex] 
 \hline
 Sinusoidal & NC-G & 7.00 & 8.50 & 15 \\ 
 & RDMSS & 35.90 & 78.50 & 41 \\ [1ex]
 \hline 
 \end{tabular}
 \label{Zakharov_deterministic_performance}
 \caption{Deterministic performance comparison for $\alpha = 0.5, \delta = 0.001$ and $\epsilon = (0.01)^5$ for compact domain $\mathcal{X} \subset \mathbb{R}^5$ aggregated over 50 individual runs per method.}
\end{table}

\section{Scalability: Effect of dimension on performance}

Our last question pertains to the robustness of the RDMSS algorithm. In particular, we devise tests to analyze the effects of dimension of the compact domain $\mathcal{X}$ on the performance of RDMSS, specifically the error between the estimated global optimum and the known solution. We expect the algorithm to lose accuracy with dimension as numerical error compounds more and more significantly. In order to determine whether or not this is the case, we test RDMSS with parameters $\epsilon = (0.01)^d$, $\alpha = 0.5$ and $\delta = 0.001$ where $d \in \{5,15,25,50\}$. We execute 50 global runs per dimension and plot the function histories on the same axes. In these experiments, we plot only a subset of the previous six objective tests, choosing only the Zakharov, Rotated Hyper-Ellipsoid and Styblinski-Tang functions as these were the ones that yielded visually notable results. \\
\indent One important consideration is that in our previous experiments, we plotted \textit{raw} function histories i.e., each function evaluation was plotted in chronological order. For the scalability tests, we plot \textit{sorted} function histories. This is simply to illustrate an overall trend that persists per dimension. Additionally, the sorted function history is a better overall representation of the global algorithm since, at the outer loop level, the global estimate is non-increasing. This was inappropriate for performance tests as illustrating the restarts was an important consideration. In the plots that follow, restarts will not be visible. 

\subsection{Zakharov}

\begin{figure}
    \centering
    \includegraphics[width=\textwidth]{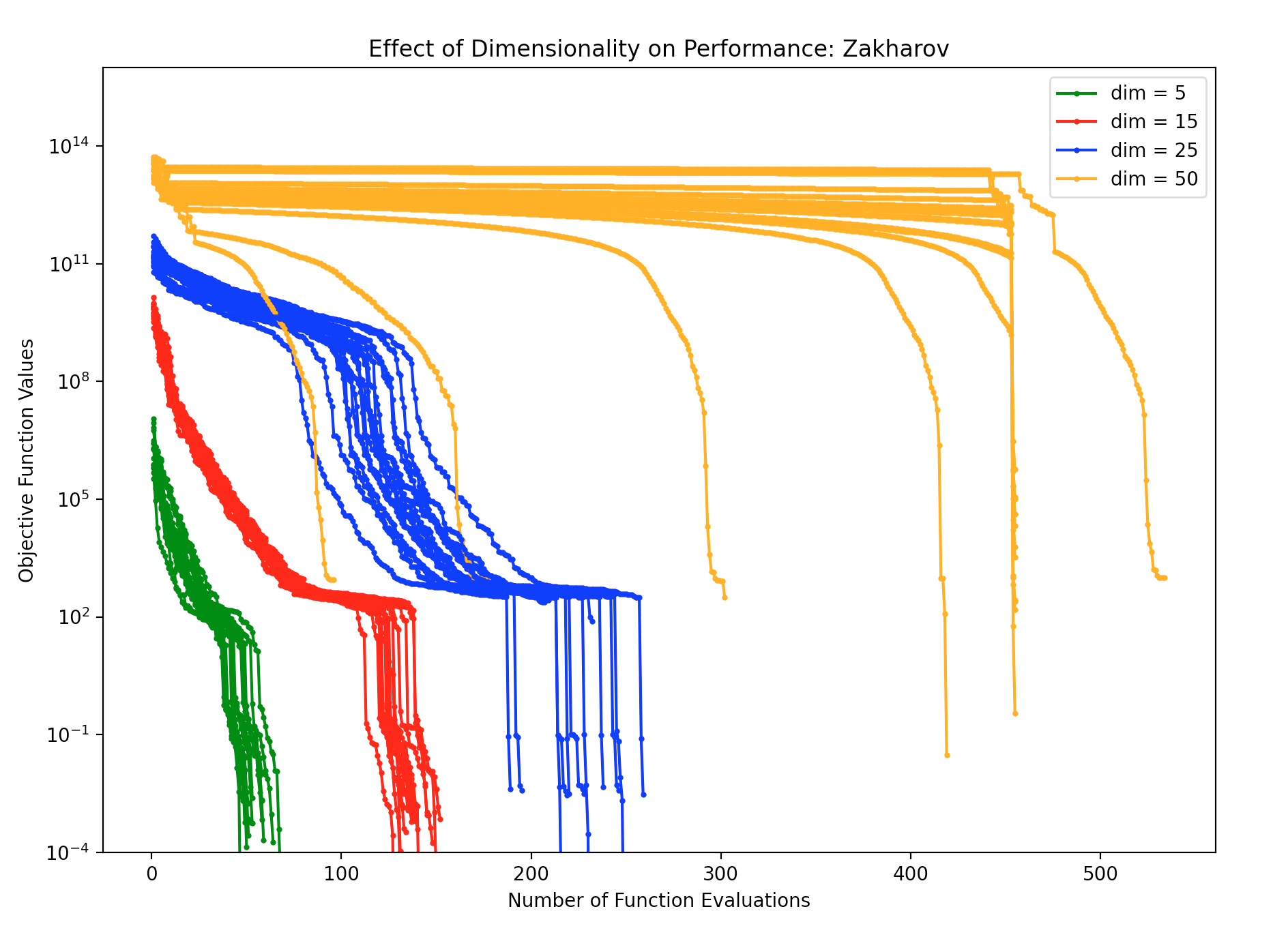}
    \caption{Sorted function histories for RDMSS in dimensions $d\in \{5,15,25,50\}$ applied to the Zakharov function}
    \label{fig:zakharov_dim}
\end{figure}

Once again we begin our experiments on the Zakharov test function. The economical use of computation that we observed for the Zakharov function during the multi-start and deterministic comparison tests begins to lead to problems in accuracy as dimension increases. What we can see from Figure \ref{fig:zakharov_dim} is that as dimension increases, the RDMSS algorithm must overcome increasingly larger ``humps". Clearly, by the form of the Zakharov function (\textit{see Appendix A}), we notice that the objective function values increase logarithmically, which means that the ``inflection error" that was alluded to in Section 4.1 compounds logarithmically as well, creating sharper bounds on the expected slope metric. This is why we can see most of the runs in higher dimensions getting stuck, with fewer and fewer ``breaking through". Additionally, for higher dimensions, even the runs that do break through their respective humps do not come close to an acceptable target region relative to their dimension. However, this may be a problem with NC-G itself. 

\subsection{Rotated Hyper-Ellipsoid}

\begin{figure}
    \centering
    \includegraphics[width=\textwidth]{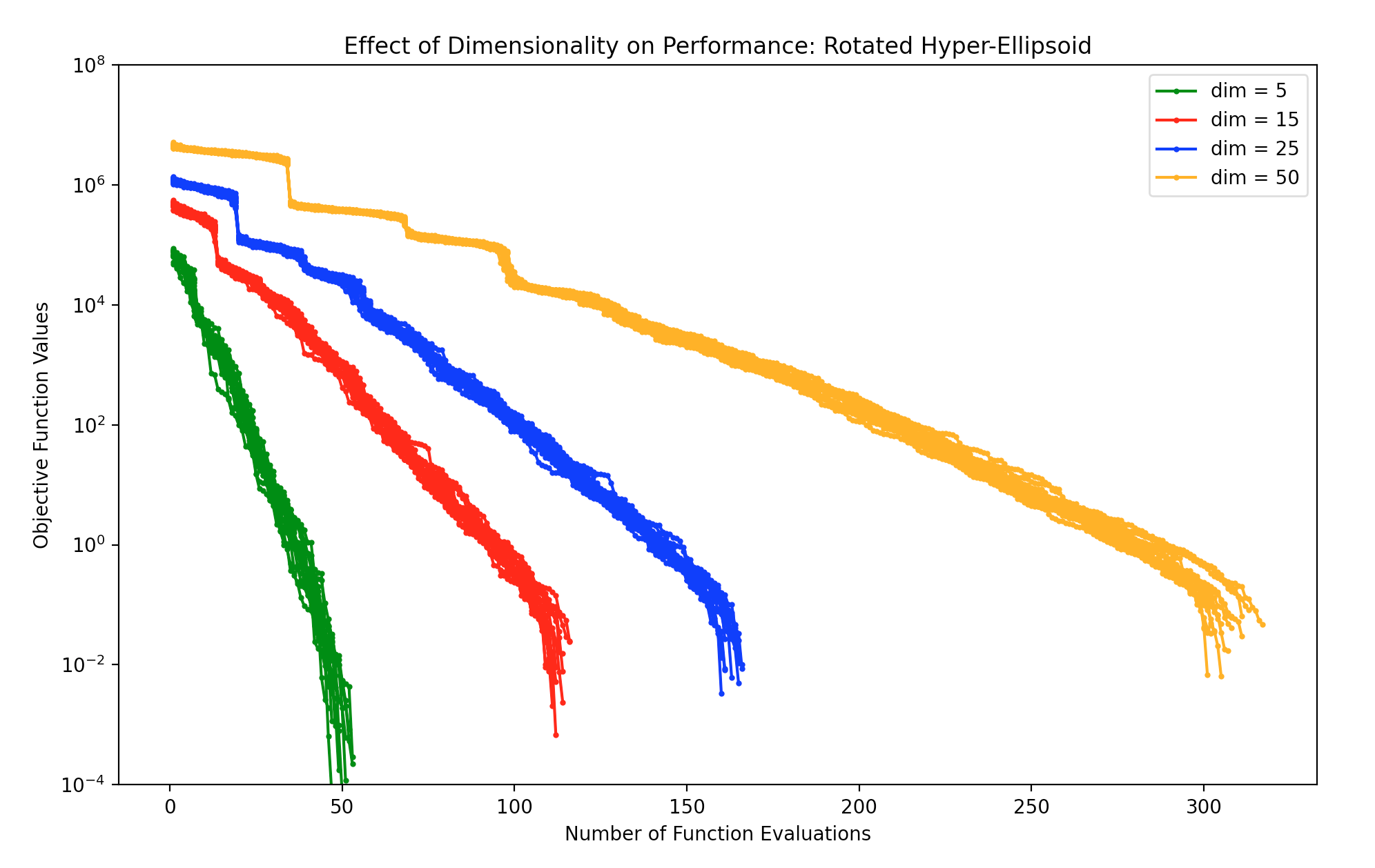}
    \caption{Sorted function histories for RDMSS in dimensions $d\in \{5,15,25,50\}$ applied to the Rotated Hyper-Ellipsoid function}
    \label{fig:rotated_hyper_ellipsoid_dim}
\end{figure}

RDMSS did not converge to the target region for Rotated Hyper-Ellipsoid function in any of the tests. However, it is of interest to visualize the linearity in error as dimension scales evident from Figure \ref{fig:rotated_hyper_ellipsoid_dim}. We can see that the sorted function histories ``cut off" at increasingly higher objective function values which suggests that the performance of RDMSS stays relatively static and we could, once again, be seeing the effects of dimensionality on NCG itself as opposed to RDMSS. 

\subsection{Styblinski-Tang}
Lastly, we observe the special case of Styblinski-Tang, which RDMSS performed well on. What is notable about the Styblinski-Tang objective is that its global optimum depends on dimension of the domain (\textit{see Appendix A}). Although it is difficult to see on the logarithmic scale, the global optimum is $-39.16599\cdot d$ where $d = \text{dim}(\mathcal{X})$. Hence, RDMSS does indeed approach the global optimum as dimension increases with an interestingly static number of function evaluations. However, we see that the accuracy begins to suffer once we go higher than dimension 15. Once again, the inherent ``impatience" of our algorithm does not care for small relative improvements and chooses to terminate earlier than required to reach the appropriate $\epsilon-$target region.

\begin{figure}
    \centering
    \includegraphics[width=\textwidth]{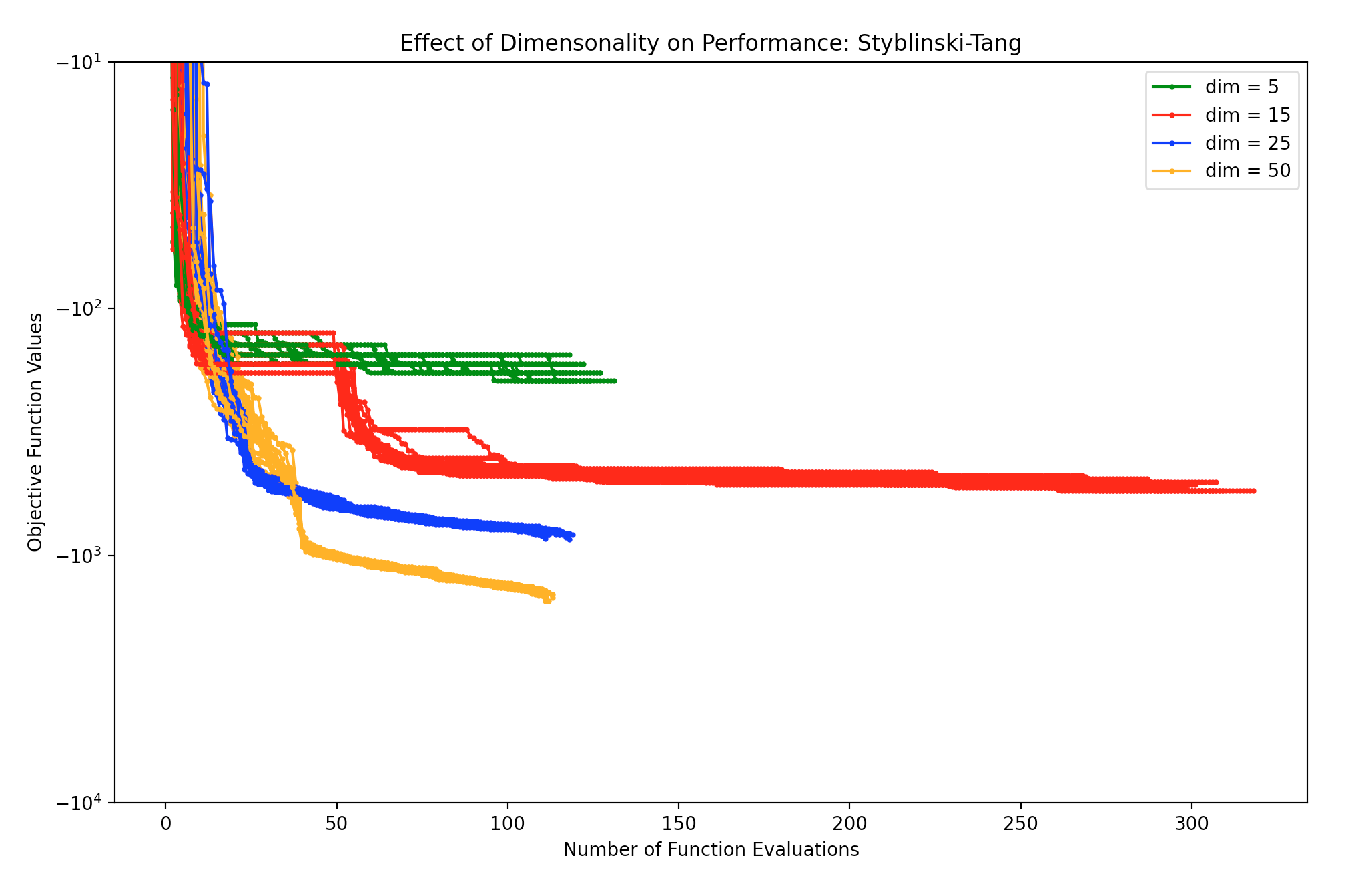}
    \caption{Sorted function histories for RDMSS in dimensions $d\in \{5,15,25,50\}$ applied to the Styblinski-Tang function}
    \label{fig:styblinski_tang_dim}
\end{figure}

\section{Conclusion}

In this thesis we consider a particular class of black-box optimization methods known as multi-start frameworks. A key question in multi-start methods is: how many inner loops (local optimization problems) to initiate and when to terminate (restart). The Dynamic Multi-start Sequential Search (DMSS) algorithm and its conceptual model, Hesitant Adaptive Search with Power-Law Improvement Distribution (HASPLID), addresses these questions but was not tested for a deterministic, gradient-based inner search. It was speculated in early observations that DMSS would perform inefficiently in such cases. \\
\indent This fact motivated the main subject of the thesis, called the Revised DMSS (RDMSS) framework. In RDMSS, we introduce a new inner loop termination criterion based on the consistent incremental improvement of local search records \textit{in addition to} the time between records that was already developed in DMSS. The new termination condition is devised by analyzing a new stochastic process called the \textit{record improvement slope process} which is based on HASPLID. We proved distributional properties of the slope process and derived a closed-form expectation. This expectation allowed us to define the metric that would become part of RDMSS' inner terminating criterion. \\
\indent Upon designing the revised framework, we subjected it to a series of experiments to answer three questions. The first was: \textit{how does RDMSS compare to DMSS for gradient-based inner searches, specifically with respect to computational expense}? We demonstrated an empirical answer to this question in the case of the Newton-Conjugate Gradient method for the inner search applied to six test functions. \\
\indent The Rosenbrock and Rotated Hyper-Ellipsoid test function experiments were indicative of a reduction in accuracy at the expense of too aggressive of a decrease in inner loop termination time. In four of the six test functions, RDMSS performed better than DMSS in terms of success rate and overall computational expense (total number of function evaluations). We found that the performance was sensitive to the shape of each specific test function. Intuitively, objective functions that experience a ``flattening" around the global optimum will \textit{trick} the inner termination criterion into terminating prematurely and thus will not allow enough time to enter the $\epsilon-$target region. The Zakharov function is a notable exception due to its ``plate-shape", meaning that the objective's entire hyper-dimensional surface is flat and thus the expected slope admits a less aggressive bound, allowing the algorithm to enter the $\epsilon-$target region.  \\
\indent The multimodal test functions, Styblinski-Tang, Shifted and Centered Sinusoidal, on the other hand, did in fact produce the predicted outcome. Namely, the Styblinski-Tang experiments showed that RDMSS was able to enter the $\epsilon-$target region 26\% of the time as opposed to DMSS' 0 \%. Additionally, it was able to do so in less overall restarts on average. The Shifted Sinusoidal function maintained its number of restarts but reduced the average number of inner search iterations and converged to the global optimum in fewer overall function evaluations while maintaining accuracy. These positive results support the hypothesis that RDMSS' more conservative slope metric was successful in reducing the overall computational expense for functions with many local minima. \\
\indent In answering our second question, we sought to compare RDMSS to NCG with no restarts to illustrate a benchmark performance for multi-start as opposed to algorithms with no restarts. We found that clearly deterministic methods are far superior in terms of success rate and computational expense for the case of unimodal objectives. RDMSS, a multi-start method, is far superior for multimodal cases. What is interesting to note is that the \textit{rates} at which both methods converge to the $\epsilon-$target region is comparable, even for unimodal cases. \\
\indent Lastly, we stress-tested RDMSS in high dimensions. We focused on the Zakharov, Rotated Hyper-Ellipsoid, and Styblinski-Tang test functions. All three tests supported the dependency of RDMSS' performance on modality of the objective test function with the scalability tests offering an extra insight about the potential drawbacks of gradient-based inner searches in general. For example, the Rotated Hyper-Ellipsoid test function was not amenable to the RDMSS algorithm and yet its error grew linearly as dimension increased. This suggests that the NC-G algorithm itself may be contributing to the weaker performance in high dimensions. \\

\indent Finally, we  close our discussion with potential future work. For further modification of the inner search in RDMSS, one can focus on the dynamic behavior of the expected slope metric as regulated by a scaling factor in the $\Tilde{p}$ approximation function for the CDF of the range measure $\rho$ in HASPLID. Modifying the algorithmic precedence of the $n_{\text{RECORD}}$ inner termination criteria could also yield interesting results. Lastly, modifying the outer search to be a more \textit{exploitative} method, such as the Gaussian Process in SOAR, could improve its performance. However, this would require a re-derivation of all the outer termination metrics as their expressions all assume independence, which would no longer hold in the case of dependent restarts.

%
%
\nocite{*}   
\bibliographystyle{plain}
\bibliography{uwthesis}

\begin{thebibliography}{10}

\bibitem{Records}
Barry~C. Arnold, N.~Balakrishnan, and H.~N. Nagaraja.
\newblock {\em Records}.
\newblock Wiley, 1998.

\bibitem{atkinson_1992}
A.~C. Atkinson.
\newblock A segmented algorithm for simulated annealing.
\newblock {\em Statistics and Computing}, 2(4):221–230, 1992.

\bibitem{grad_trapped}
Hemayet~Ahmed Chowdhury, Md. Azizul~Haque Imon, Anisur Rahman, Aisha Khatun,
  and Md.~Saiful Islam.
\newblock A continuous space neural language model for bengali language.
\newblock {\em CoRR}, abs/2001.05315, 2020.

\bibitem{durrett_2010}
Richard Durrett.
\newblock {\em Probability: theory and examples}.
\newblock Cambridge Univ. Press, 2010.

\bibitem{gildewell_ng_hensel}
M.e. Gildewell, K.t. Ng, and E.~Hensel.
\newblock A combinatorial optimization approach as a pre-processor for
  impedance tomography.
\newblock {\em Proceedings of the Annual International Conference of the IEEE
  Engineering in Medicine and Biology Society Volume 13: 1991}.

\bibitem{gupta_1984}
Ramesh~C. Gupta.
\newblock Relationships between order statistics and record values and some
  characterization results.
\newblock {\em Journal of Applied Probability}, 21(2):425–430, 1984.

\bibitem{harmelen_lifschitz_porter_2010}
Frank~Van Harmelen, Vladimir Lifschitz, and Bruce Porter.
\newblock {\em Handbook of knowledge representation}.
\newblock Elsevier, 2010.

\bibitem{Horst1995HandbookOG}
Reiner Horst, Panos~M. Pardalos, and H.~Edwin Romeijn.
\newblock {\em Handbook of Global Optimization}.
\newblock Springer, 1995.

\bibitem{Horst2000IntroductionTG}
Reiner Horst, Panos~M. Pardalos, and Nguyen~V. Thoai.
\newblock {\em Introduction to Global Optimization}.
\newblock Springer, US, 2000.

\bibitem{annealing}
S.~Kirkpatrick, C.~D. Gelatt, and M.~P. Vecchi.
\newblock Optimization by simulated annealing.
\newblock {\em Science}, 220(4598):671--680, 1983.

\bibitem{li_lim}
H.~Li and A.~Lim.
\newblock A metaheuristic for the pickup and delivery problem with time
  windows.
\newblock In {\em Proceedings 13th IEEE International Conference on Tools with
  Artificial Intelligence. ICTAI 2001}, pages 160--167, 2001.

\bibitem{glob_opt_based_local_search}
Marco Locatelli and Fabio Schoen.
\newblock {Global optimization based on local searches}.
\newblock {\em Annals of Operations Research}, 240(1):251--270, May 2016.

\bibitem{Marti2003}
Rafael Mart{\'i}.
\newblock {\em Multi-Start Methods}.
\newblock Springer US, Boston, MA, 2003.

\bibitem{SOAR}
L.~Mathesen, G.~Pedrielli, S.H. Ng, and Z.B. Zabinsky.
\newblock Stochastic optimization with adaptive restart: a framework for
  integrated local and global learning.
\newblock {\em Journal of Global Optimization}, 79:87--110, 2021.

\bibitem{nagaraja_nevzorov_1997}
H.N. Nagaraja and V.B. Nevzorov.
\newblock On characterizations based on record values and order statistics.
\newblock {\em Journal of Statistical Planning and Inference}, 63(2):271–284,
  1997.

\bibitem{neumann_witt_2007}
Frank Neumann and Carsten Witt.
\newblock Runtime analysis of a simple ant colony optimization algorithm.
\newblock {\em Algorithmica}, 54(2):243–255, 2007.

\bibitem{StochProc}
Sheldon~M. Ross.
\newblock {\em Stochastic Processes}.
\newblock Wiley India, 2016.

\bibitem{shifted_sine_paper}
Yanfang Shen, Seksan Kiatsupaibul, Zelda Zabinsky, and Robert Smith.
\newblock An analytically derived cooling schedule for simulated annealing.
\newblock {\em Journal of Global Optimization}, 38:333--365, 06 2007.

\bibitem{local_opt_review}
Gerhard Venter.
\newblock {\em Review of Optimization Techniques}.
\newblock Wiley, 12 2010.

\bibitem{GlobalOpt}
Zelda~B. Zabinsky.
\newblock {\em Stochastic Adaptive Search for Global Optimization.}
\newblock Kluwer Academic Publishers, 2003.

\bibitem{HASPLID}
Zelda~B. Zabinsky, David Bulger, and Charoenchai Khompatraporn.
\newblock Stopping and restarting strategy for stochastic sequential search in
  global optimization.
\newblock {\em J. Global Optimization}, 46:273--286, 02 2010.

\end{thebibliography}
%
%
\appendix
\raggedbottom\sloppy
 

\section{Objective Test Functions}
Note that $d$ denotes the dimension of $\mathcal{X}$, the domain of the respective test function $f$. 
\begin{enumerate}
    \item \textit{Zakharov Function}
    \begin{equation}
    f(x) = \sum_{i=1}^d x_i^2 + \left( \sum_{i=1}^d 0.5ix_i \right)^2 + \left( \sum_{i=1}^d 0.5ix_i \right)^4
    \end{equation}
    with domain $\mathcal{X} = [-5,10]^d$ and $f(x^*) = 0$ for $x^* = (0,0, \cdots, 0)$.
    \item \textit{Rosenbrock Function}
    \begin{equation}
    f(x) = \sum_{i=1}^d\left( 100(x_{i+1} - x_i^2)^2 + (x_i - 1)^2 \right)
    \end{equation}
    with domain $\mathcal{X} = [-2.048, 2.048]^d$ and $f(x^*) = 0$ for $x^* = (1,1,\cdots, 1)$.
    \item \textit{Rotated Hyper-Ellipsoid}
    \begin{equation}
    f(x) = \sum_{i=1}^d\sum_{j=1}^i x_j^2
    \end{equation}
    with domain $\mathcal{X} = [-65.536,65.536]^d$ with global minimum $f(x^*) = 0$ for $x^* = (0,0,\cdots,0)$. 
    \item \textit{Styblinski-Tang}
    \begin{equation}
        f(x) = \frac{1}{2}\sum_{i=1}^5 (x_i^4 - 16x_i^2 + 5x_i)
    \end{equation}
    with domain $\mathcal{X} = [-5,5]^d$ and global minimum $f(x^*) = -39.16599 \cdot d $ at $x^* = (-2.903534, \cdots, -2.903534)$.
    \item \textit{Shifted Sinusoidal}
    \begin{equation}
        f(x) = -2.5\prod_{i=1}^d\sin(x_i + 60) - \prod_{i=1}^d\sin(5(x_i+60))
    \end{equation}
    with domain $\mathcal{X} = [-90,90]^d$ and global minimum $f(x^*) = -3.5 $ at $x^* = (30,30, \cdots, 30)$. 
    \item \textit{Centered Sinusoidal}
    \begin{equation}
        f(x) = -2.5\prod_{i=1}^d\sin(x_i + 90) - \prod_{i=1}^d\sin(5(x_i+90))
    \end{equation}
    with domain $\mathcal{X} = [-90,90]^d$ and global minimum $f(x^*) = -3.5 $ at $x^* = (0,0, \cdots, 0)$. 
    
\end{enumerate}
\comment{
\section{Hesitant Adaptive Search}
 
\medskip

Hesitant adaptive search (HAS) extends the idea of pure adaptive search (PAS \texit{c.f.} \cite{GlobalOpt}) by introducing a bettering probability. This bettering probability $b(y)$ dictates when the search algorithm will generate a point in the improving region or ``hesitate" i.e. stay at its current location in the effective domain. 

\medskip 

Returning to the problem (1.1), we will measure performance of our minimization by computing the \textit{number of iterations to first sample an objective function value of $y$ or less}. Since we have $\mu$ as the measure on our sample space, we begin the algorithm by selecting $X_0$, the first sample point, by $\mu$. For every subsequent iteration, the algorithm hesitates at the current objective function value $y$ with probability $1 - b(y)$ or improves with probability $b(y)$ (which we assume is measurable and bounded away from zero). The steps of HAS are as follows: 
\begin{quote}
\textbf{Hesitant Adaptive Search$(\mu, b(y))$} \textit{c.f.} \cite{GlobalOpt}:
    \begin{enumerate}
        \item[\textbf{Step 0.}] Initialize $X_0 \in S$ according to probability measure $\mu$ on measurable space $(S,\mathcal{F})$. Set $k=0$ and $Y_0 := f(X_0)$. 
        \item[\textbf{Step 1.}] Generate $X_{k+1}$ according to the normalized restriction of $\mu$ on the level set $S_k = \{x: x \in S, f(x) < Y_k\}$ with probability $b(Y_k)$; otherwise set $X_{k+1} = X_k$. Set $Y_{k+1} := f(X_{k+1})$. 
        \item[\textbf{Step 2.}] If a stopping criterion is met, stop; otherwise, increment $k$ and return to \textbf{Step 1.}
    \end{enumerate}
\end{quote}
To introduce some further notation: $N(y)$ denotes the number of iterations required to first achieve an objective function value of $y \in \mathbb{R}$ or less. Alternatively, we can think of $N(y)$ as the number of iterations until HAS lands in the termination region $T_y$.

\medskip

The distribution of $N(y)$ can be analytically computed by considering $M(y) = N(y) - 1$.  To ensure that $\mathbb{E}[N(y)] < \infty$, we assume that $\rho(T_y) > 0 ~ \forall x \in \{x \in T_y: x < y\}$. Finally, we return back to our definition of $p$, the CDF of $\rho$ where $$p(y) = \rho((-\infty,y]) = \mu(S(y)) \textrm{ and } S(y) = \{x \in S: f(x) \leq y, ~ y \in \mathbb{R} \}.$$

\section{Continuous HAS Distribution}
 
We let $Y_{R(k)}$ denote the $k$th HAS record, i.e. the $k$th improving objective function value. If it takes $N(y)$ distinct objective function values in a run of HAS before the termination level $y$ is reached, then the records of the run are given by: $Y_{R(1)}, Y_{R(2)},\hdots, Y_{R(N(y))}$. \\ 

\medskip 

\noindent \textbf{Lemma 3.1}(\textit{c.f} \cite{GlobalOpt}) ~ The stochastic process $\{Y_{R(k)}: k \in \mathbb{N}\}$ of HAS records is equal in distribution to the stochastic process $\{W_k, k\in \mathbb{N}\}$ of PAS, i.e. $$\{W_k\} \sim \{Y_{R(k)}\} $$ when using the same range probability measure $\rho$. \\ 
\textit{Remark}. A key step in the proof is noticing that for $x,y \in \mathbb{R},$ such that $y \leq x$, $$\mathbb{P}[Y_{R(k+1)} \leq y|Y_{R(k)} = x] = \frac{p(y)}{p(x)} = \mathbb{P}[W_{k+1} \leq y|W_k = x].$$ 

We return again to the number of iterations to convergence \textit{just prior} to achieving the termination region $[-\infty, y),~ M(y) = N(y) - 1$. If $M_i$ is the number of hesitations at level $i$ prior to a new record, then the total number of HAS iterations prior to reaching level $y$ is: $$M(y) = \sum_{i=1}^{T_y} M_i $$ where $M_i \sim \text{Geo}(b(Y_{R(i)})).$ \\

\noindent\textit{Note:} Theorem 2.4 (\textit{c.f.} \cite{GlobalOpt}) states that the process formed by $(M(z))_{z \geq 0}$ for PAS is a non-homogeneous Poisson process with mean value function $$m(z) = \log\left(p\left(\frac{y^* + zy_*}{1+z}\right) \right)^{-1}. $$ Intuitively this makes sense that introducing hesitations into PAS makes HAS in the continuous case a \textit{marked Poisson Process}.\\

Now we list two key results regarding the quantity $M(y)$: 
\begin{enumerate}
    \item[(1)](\textit{c.f.} \cite{GlobalOpt} Thm 3.2) $$\mathbb{E}\left[ z^{M(y)} \right] = \exp{\int_y^\infty \frac{z-1}{(zb(t) - (z-1))p(t)} d\rho(t) } $$ 
    \item[(2)](\textit{c.f.} \cite{GlobalOpt} Cor. 3.3) $$\mathbb{E}[N(y)] = 1 + \mathbb{E}[M(y)] = 1 + \int_y^\infty \frac{d\rho(t)}{b(t)p(t)} $$ $$\text{Var}[N(y)] = \text{Var}[M(y)] = \int_y^\infty \left(\frac{2}{b(t)} - 1\right) \frac{2\rho(t)}{b(t)p(t)} $$
\end{enumerate}
 }
 
\section{Record Value Theory}

A standard record value process is given by $(X_k)_{k\in\mathbb{N}}$ which are i.i.d. to $X$ with CDF $F$. An observation $X_j$ is called a \textit{lower-record value} or \textit{record} if $X_j < X_i$ for all $i < j$. Note that this is equivalent to the definition of an \textit{order statistic}. Assuming discrete-time, our index set $\{k \in \mathbb{N}\}$ denotes the chronology in which our observations $(X_k)$ appear. \\
\indent Thus, the \textit{record time sequence}, $(R(k))_{k \in \mathbb{N}}$ is defined by: $$R(0) = 1 ~ ~ ~ ~ ~ ~ ~ ~ \text{w.p.  } 1 $$ and $$ R(k) = \min\{j: X_j < X_{R(k)-1}\}$$ for $k \geq 1$. This is simply the collection of indices upon which records appear. Note that each $R(k)$ is itself a random variable. \textit{Note:} In general, it is assumed that the process $(X_k)$ does not admit an \textit{unbreakable} record i.e., $X_k$ does not have a lower bound. However, in the context of optimization, existence of an optimum clearly indicates admittance of an unbreakable record, namely, the optimal value itself. \\
\indent The \textit{record incrememt process} $\{Y_{R(k)} - Y_{R(k-1)}, k \geq 1\}$, sometimes called a jump process \cite{Records}, is defined as above where $Y_{R(1)} - Y_{R(0)} = Y_{R(1)}$ since we initialize the $0$th record at $\infty$. Of similar importance is the the \textit{inter-record time} sequence $\{R(k) - R(k-1), k \geq 1\}$. If $(X_k)$ are i.i.d \textit{continuous} random variables, \cite{Records} names this setting the \textit{classical record model}. 

\section{Classical Record Model}

What makes the classical record model work is the assumption that observations are exponential. If $X = (X_j) \sim $Exp$(1)$ random variables then the lack of memory property yields $\{Y_{R(k)} - Y_{R(k-1)}, k \geq 1\} \sim $Exp$(1)$ and thus, \begin{equation}
    Y_{R(k)} \sim G(k+1,1).
\end{equation}
Further, if $X$ has a \textit{continuous} CDF $F$, then 
\begin{equation}
    H(X) := -\log(1 - F(X))
\end{equation} is the distribution of the standard exponential random variable. Now, since the sequence of records is monotone, we can obtain the following expression: 
\begin{equation}
    \mathbb{P}(Y_{R(k)} > r^*) = e^{-r^*} \cdot \sum_{k=0}^n(r^*)^k/k!, ~ ~ ~ r^* > 0.
\end{equation}

\end{document}